\documentclass[12pt,reqno]{amsart}
\usepackage{amsmath,amssymb,amsthm}
\usepackage{geometry}
\usepackage{graphicx}
\usepackage{cite}
\usepackage{hyperref}

\newtheorem{theorem}{Theorem}[section]
\newtheorem{lemma}[theorem]{Lemma}
\newtheorem{proposition}[theorem]{Proposition}
\newtheorem{corollary}[theorem]{Corollary}

\theoremstyle{definition}
\newtheorem{definition}[theorem]{Definition}
\newtheorem{example}[theorem]{Example}
\newtheorem{remark}[theorem]{Remark}

\DeclareMathOperator{\Lip}{Lip}
\DeclareMathOperator{\diam}{diam}

\begin{document}

\title[Contraction Divisors and Coupled Self-Affine Sets]{Contraction Divisors and the Hausdorff Dimension of Coupled Self-Affine Sets}

\author{Amal P. S.}
\address{APJ Abdul Kalam Technological University, Thiruvananthapuram, Kerala, India; Rajagiri School of Engineering and Technology, Kochi, Kerala, India}
\email{amalb9111@gmail.com}

\author{Ramkumar P. B.}
\address{APJ Abdul Kalam Technological University, Thiruvananthapuram, Kerala, India; Rajagiri School of Engineering and Technology, Kochi, Kerala, India}
\email{ramkumar\_pb@rajagiritech.edu.in}

\author{Vinod Kumar P. B.}
\address{APJ Abdul Kalam Technological University, Thiruvananthapuram, Kerala, India; Muthoot Institute of Technology and Science, Kochi, Kerala, India}
\email{vinodkumarpb@mgits.ac.in}

\subjclass[2020]{Primary 28A80, 28A78; Secondary 37C45, 54E40}
\keywords{Contraction divisors, Iterated function systems, Self-affine sets, Hausdorff dimension, Open set condition}
\thanks{Amal P. S. is the corresponding author.}

\begin{abstract}
Motivated by zero divisors in ring theory, we introduce \emph{contraction divisors}: non-contractive mappings $f$ on a metric space for which there exists a non-contractive mapping $g$ such that $f \circ g$ and $g \circ f$ are contractions. We characterize affine contraction divisors on $\mathbb{R}^m$ via singular values and establish fixed-point and convergence properties on complete metric spaces. We then introduce coupled iterated function systems, in which two families of contraction divisors act alternately on a pair of compact sets alongside classical contractions. For a class of coupled self-affine systems on $\mathbb{R}^m$ combining affine contraction divisors with compatible similarities, we prove the existence and uniqueness of coupled attractors, derive an exact Hausdorff dimension formula under a coupled open set condition, and show that the Hausdorff, packing, and box-counting dimensions coincide.
\end{abstract}

\maketitle

\section*{Introduction}

On a complete metric space $(X, d)$, Banach's Fixed Point Theorem \cite{Banach1922} guarantees a unique fixed point for every contraction $f: X \to X$. Building on this principle, Hutchinson \cite{Hutchinson1981} proved that any finite collection of contractions $\{f_1, \dots, f_k\}$ admits a unique non-empty compact invariant set $A = \bigcup_{i=1}^k f_i(A)$, the attractor of the corresponding Iterated Function System (IFS) \cite{BarnsleyDemko1985,Barnsley2014,Falconer2014}. When an IFS of similarities on $\mathbb{R}^m$ with contraction ratios $r_1, \dots, r_k$ satisfies the Open Set Condition (OSC), the Hausdorff dimension of its attractor $A$ equals the similarity dimension $s$ defined by the Moran equation $\sum_{i=1}^k r_i^s = 1$ \cite{Moran1946}; moreover, Schief \cite{Schief1994} proved that the OSC is equivalent to $\mathcal{H}^s(A) > 0$. Mauldin and Williams \cite{MauldinWilliams1988} extended this theory to graph-directed constructions of similarities in $\mathbb{R}^m$, determining the Hausdorff dimension via the spectral radius of a weighted incidence matrix, while Nussbaum, Priyadarshi, and Verduyn Lunel \cite{Nussbaum2012} used positive linear operators to characterize the Hausdorff dimension of invariant sets for systems of infinitesimal similitudes.

For self-affine sets, where maps contract different directions at different rates, no single scalar equation determines the dimension in general. Falconer \cite{Falconer1988} introduced the affinity dimension as an upper bound and showed that, under a norm bound on the linear parts, the Hausdorff dimension equals the minimum of the affinity dimension and the ambient dimension $m$ for almost all translation vectors. Subsequent work established explicit conditions under which this equality holds, beginning with totally disconnected planar sets \cite{HueterLalley1995} and advancing through Furstenberg measures, projections, and positive subsystems under separation, irreducibility, or generic linear-part hypotheses \cite{FalconerKempton2018,Barany2018,BaranyHochmanRapaport2019,MorrisShmerkin2019}, as well as for overlapping planar systems \cite{BaranyRamsSimon2016,HochmanRapaport2022}. More recently, Feng \cite{Feng2023} proved exact dimensionality and a Ledrappier--Young formula for ergodic invariant measures of contracting---and more generally, average contracting---affine IFSs on $\mathbb{R}^m$, while Rapaport \cite{Rapaport2026} established Falconer's formula for diagonal self-affine sets on $\mathbb{R}^m$ under coordinate-wise exponential separation and a non-conformality condition, using non-conformal partitions.

Alongside these developments, the classical contraction condition has been generalized in several directions, including Kannan, Edelstein, and Rakotch contractions \cite{Kannan1969,Edelstein1962,Rakotch1962,Rhoades1977}. In the IFS setting, Le\'sniak, Snigireva, Strobin, and Vince \cite{LesniakSnigirevaStrobinVince2025} recently showed that even IFSs on $\mathbb{R}^m$ whose maps have Lipschitz constant strictly greater than $1$ under every equivalent metric can admit a unique attractor. Another line of work concerns mappings with contractive higher iterates, whether global \cite{Bryant1968} or point-dependent \cite{Sehgal1969}, and their extension to IFSs \cite{AmalGIFS2026}. In \cite{AmalSelfAffine2026}, exact Hausdorff dimension formulas were obtained under the OSC for self-affine systems combining affine maps whose length-$n$ compositions are similarity contractions with similarities whose linear parts commute with $M^\top M$ for each affine linear part $M$.

If we view the collection of contractions $\mathcal{C}(X)$ on $(X, d)$ as an analogue of zero inside the monoid of self-mappings under composition, then a non-contractive mapping $f$ with a contractive iterate ($f^n \in \mathcal{C}(X)$ for some $n \ge 2$) corresponds to a nilpotent element. Motivated by zero divisors in ring theory, we consider non-contractive mappings $f$ and $g$ whose individual iterates $f^n$ and $g^n$ need not be contractions for any $n \ge 1$, yet whose compositions $f \circ g$ and $g \circ f$ lie in $\mathcal{C}(X)$. We call such mappings \emph{contraction divisors}.

The paper is organized as follows. In Section~\ref{sec:cd}, we define contraction divisors, characterize affine contraction divisors on $\mathbb{R}^m$ via singular values, and prove a fixed-point and convergence result for alternating iterations on complete metric spaces. In Section~\ref{sec:cd_ifs}, we introduce coupled contraction-divisor iterated function systems on complete metric spaces, combining alternating families of contraction divisors with classical contractions; we show that a unique coupled attractor exists whenever at least one family of classical contractions is empty, whereas a coupled attractor need not exist in general. In Section~\ref{sec:dim}, we specialize to coupled self-affine systems on $\mathbb{R}^m$ combining compatible similarity contractions with affine contraction divisors whose pairwise compositions are similarity contractions. We establish the existence and uniqueness of their coupled attractors, derive an exact formula for their Hausdorff dimension $s$ under a coupled open set condition, and show that their Hausdorff, packing, and box-counting dimensions coincide; moreover, whenever $s > 0$, the attractors support Ahlfors $s$-regular measures.

\section{Contraction Divisors in Metric Spaces}\label{sec:cd}

Let $(X, d)$ be a metric space. For any mapping $f: X \to X$, its Lipschitz constant is defined by
\[
\Lip(f) = \sup_{\substack{x, y \in X \\ x \neq y}} \frac{d(f(x), f(y))}{d(x, y)}.
\]
We denote by $\mathcal{C}(X)$ the collection of all classical contractions on $X$, that is, mappings $f: X \to X$ satisfying $\Lip(f) < 1$. Throughout this paper, $\mathbb{N} = \{1, 2, \dots\}$, $\mathbb{R}^m$ is equipped with the standard Euclidean metric $d(x, y) = \|x - y\|$, and $\|M\|$ denotes the induced Euclidean operator norm of a real matrix $M$.

\begin{definition}\label{def:cd}
Let $(X, d)$ be a metric space and let $f: X \to X$ be a mapping such that $f \notin \mathcal{C}(X)$.
\begin{enumerate}
    \item $f$ is called a \emph{left contraction divisor} if there exists $g: X \to X$ with $g \notin \mathcal{C}(X)$ such that $f \circ g \in \mathcal{C}(X)$.
    \item $f$ is called a \emph{right contraction divisor} if there exists $g: X \to X$ with $g \notin \mathcal{C}(X)$ such that $g \circ f \in \mathcal{C}(X)$.
    \item $f$ is called a \emph{contraction divisor} if there exists $g: X \to X$ with $g \notin \mathcal{C}(X)$ such that both $f \circ g \in \mathcal{C}(X)$ and $g \circ f \in \mathcal{C}(X)$. In this case, $(f, g)$ is called a \emph{contraction divisor pair}; it is called \emph{proper} if $f^n \notin \mathcal{C}(X)$ and $g^n \notin \mathcal{C}(X)$ for all $n \in \mathbb{N}$.
\end{enumerate}
\end{definition}

If $f \notin \mathcal{C}(X)$ satisfies $f^2 \in \mathcal{C}(X)$, then $(f, f)$ is a contraction divisor pair that is not proper. The following examples illustrate proper contraction divisor pairs.

\begin{example}\label{ex:basic_proper}
Define $f, g: \mathbb{R}^2 \to \mathbb{R}^2$ by
\[
f(x, y) = \left(2x, \, \frac{y}{4}\right), \qquad g(x, y) = \left(\frac{x}{4}, \, 2y\right).
\]
For every $n \in \mathbb{N}$, $d(f^n(1, 0), f^n(0, 0)) = 2^n \ge 2$ and $d(g^n(0, 1), g^n(0, 0)) = 2^n \ge 2$ (in particular, both sequences tend to $\infty$ as $n \to \infty$). Since $d((1, 0), (0, 0)) = d((0, 1), (0, 0)) = 1$, it follows that $f^n, g^n \notin \mathcal{C}(\mathbb{R}^2)$ for all $n \in \mathbb{N}$. However,
\[
(f \circ g)(x, y) = (g \circ f)(x, y) = \left(\frac{x}{2}, \, \frac{y}{2}\right),
\]
which is a contraction on $\mathbb{R}^2$ with $\Lip(f \circ g) = \Lip(g \circ f) = 1/2$. Hence $(f, g)$ is a proper contraction divisor pair.
\end{example}

While a mapping $f$ on a complete metric space with $f^n \in \mathcal{C}(X)$ for some $n \in \mathbb{N}$ always admits a unique fixed point \cite{Bryant1968}, the individual mappings in a proper contraction divisor pair need not possess any fixed points.

\begin{example}\label{ex:no_fixed_points}
Define $f, g: \mathbb{R}^2 \to \mathbb{R}^2$ by
\[
f(x, y) = \left(x + 1, \, \frac{y}{3}\right), \qquad g(x, y) = \left(\frac{x}{3}, \, y + 1\right).
\]
The equation $f(x, y) = (x, y)$ requires $x + 1 = x$, which has no solution; similarly, $g(x, y) = (x, y)$ has no solution. Thus, neither $f$ nor $g$ admits a fixed point in $\mathbb{R}^2$. Moreover, $d(f^n(0, 0), f^n(1, 0)) = 1$ and $d(g^n(0, 0), g^n(0, 1)) = 1$ for all $n \in \mathbb{N}$, so $f^n, g^n \notin \mathcal{C}(\mathbb{R}^2)$ for all $n \in \mathbb{N}$. However,
\[
(f \circ g)(x, y) = \left(\frac{x}{3} + 1, \, \frac{y}{3} + \frac{1}{3}\right), \qquad (g \circ f)(x, y) = \left(\frac{x}{3} + \frac{1}{3}, \, \frac{y}{3} + 1\right)
\]
are contractions on $\mathbb{R}^2$ with ratio $1/3$, possessing the unique fixed points $x_0 = (3/2, 1/2)$ and $y_0 = (1/2, 3/2)$, respectively. Thus $(f, g)$ is a proper contraction divisor pair.
\end{example}

\begin{remark}\label{rem:non_remetrizable}
Proper contraction divisors cannot, in general, be realized as classical contractions by passing to a topologically equivalent complete metric on $X$. Indeed, a mapping admitting an unbounded orbit (Example~\ref{ex:basic_proper}) or no fixed point (Example~\ref{ex:no_fixed_points}) cannot be a contraction under any topologically equivalent complete metric on $X$, since a classical contraction on a complete metric space must possess a unique globally attracting fixed point.
\end{remark}

We next establish that uniformly expanding mappings cannot be contraction divisors.

\begin{proposition}\label{prop:obstruction}
Let $(X, d)$ be a metric space containing at least two points, and let $f: X \to X$ be a mapping satisfying
\begin{equation}\label{eq:expanding}
d(f(x), f(y)) \ge m \, d(x, y) \quad \text{for all } x, y \in X,
\end{equation}
for some constant $m \ge 1$.
\begin{enumerate}
    \item The mapping $f$ cannot be a left contraction divisor.
    \item If $f$ is surjective, then $f$ cannot be a right contraction divisor; more generally, if $f(X)$ is dense in $X$, then no continuous mapping $g \notin \mathcal{C}(X)$ can satisfy $g \circ f \in \mathcal{C}(X)$.
\end{enumerate}
\end{proposition}

\begin{proof}
(1) Suppose, for the sake of contradiction, that there exists $g: X \to X$ with $g \notin \mathcal{C}(X)$ such that $f \circ g \in \mathcal{C}(X)$ with contraction constant $k \in [0, 1)$. Applying \eqref{eq:expanding} to the points $g(x)$ and $g(y)$ gives
\[
m \, d(g(x), g(y)) \le d(f(g(x)), f(g(y))) \le k \, d(x, y) \quad \text{for all } x, y \in X.
\]
Dividing by $m \ge 1$ yields $d(g(x), g(y)) \le (k/m) d(x, y)$. Since $k/m \le k < 1$, we obtain $g \in \mathcal{C}(X)$, a contradiction.

(2) Suppose $g \circ f \in \mathcal{C}(X)$ with contraction constant $k \in [0, 1)$. Given any $u, v \in f(X)$, choose $x, y \in X$ such that $u = f(x)$ and $v = f(y)$. By \eqref{eq:expanding}, $d(x, y) \le \frac{1}{m} d(f(x), f(y)) = \frac{1}{m} d(u, v)$, and therefore
\[
d(g(u), g(v)) = d(g(f(x)), g(f(y))) \le k \, d(x, y) \le \frac{k}{m} d(u, v).
\]
If $f$ is surjective, then $f(X) = X$, so $\Lip(g) \le k/m < 1$, contradicting $g \notin \mathcal{C}(X)$. If $f(X)$ is dense in $X$ and $g$ is continuous, the inequality $d(g(u), g(v)) \le (k/m) d(u, v)$ extends by continuity from $f(X)$ to all $u, v \in X$, again yielding $g \in \mathcal{C}(X)$, a contradiction.
\end{proof}

For affine mappings on $\mathbb{R}^m$, Proposition~\ref{prop:obstruction} leads to a singular-value characterization.

\begin{theorem}\label{thm:affine_cd_char}
Let $f: \mathbb{R}^m \to \mathbb{R}^m$ ($m \ge 1$) be an affine mapping given by $f(x) = Ax + a$, and let $\sigma_{\min}(A)$ and $\sigma_{\max}(A)$ denote the smallest and largest singular values of $A$, respectively. The following statements are equivalent:
\begin{enumerate}
    \item $f$ is a contraction divisor;
    \item $f$ is a left contraction divisor;
    \item $f$ is a right contraction divisor;
    \item $\sigma_{\min}(A) < 1 \le \sigma_{\max}(A)$.
\end{enumerate}
\end{theorem}

\begin{proof}
The implications $(1) \Rightarrow (2)$ and $(1) \Rightarrow (3)$ are immediate from Definition~\ref{def:cd}.

To show $(2) \Rightarrow (4)$ and $(3) \Rightarrow (4)$, suppose $f$ is either a left or a right contraction divisor. Then $f \notin \mathcal{C}(\mathbb{R}^m)$, so $\sigma_{\max}(A) = \Lip(f) \ge 1$. If $\sigma_{\min}(A) \ge 1$, then $A$ is invertible (so $f$ is surjective) and $\|f(x) - f(y)\| = \|A(x - y)\| \ge \sigma_{\min}(A)\|x - y\| \ge \|x - y\|$ for all $x, y \in \mathbb{R}^m$. By Proposition~\ref{prop:obstruction}, $f$ can be neither a left nor a right contraction divisor, a contradiction. Hence $\sigma_{\min}(A) < 1 \le \sigma_{\max}(A)$.

To prove $(4) \Rightarrow (1)$, assume $\sigma_{\min}(A) < 1 \le \sigma_{\max}(A)$, so $f \notin \mathcal{C}(\mathbb{R}^m)$.
\begin{itemize}
    \item If $\sigma_{\min}(A) > 0$, then $A$ is invertible. Choose $c \in [\sigma_{\min}(A), 1)$ and define $g(x) = cA^{-1}x$. Then $\Lip(g) = \|cA^{-1}\| = c/\sigma_{\min}(A) \ge 1$, so $g \notin \mathcal{C}(\mathbb{R}^m)$, while $(f \circ g)(x) = cx + a$ and $(g \circ f)(x) = cx + cA^{-1}a$ both have Lipschitz constant $c < 1$.
    \item If $\sigma_{\min}(A) = 0$, then $A$ is singular, so $\ker(A) \neq \{0\}$ and $(\operatorname{Im}(A))^\perp \neq \{0\}$. Choose unit vectors $v \in \ker(A)$ and $u \in (\operatorname{Im}(A))^\perp$, and define $g(x) = \langle x, u \rangle v$. Then $\|g(x) - g(y)\| = |\langle x - y, u \rangle| \le \|x - y\|$ with equality when $x - y = u$, so $\Lip(g) = 1$ and $g \notin \mathcal{C}(\mathbb{R}^m)$. Meanwhile, $(f \circ g)(x) = \langle x, u \rangle Av + a = a$ and $(g \circ f)(x) = \langle Ax + a, u \rangle v = \langle a, u \rangle v$ are constant maps, hence contractions with Lipschitz constant $0$.
\end{itemize}
In both cases, $(f, g)$ is a contraction divisor pair.
\end{proof}

In particular, Theorem~\ref{thm:affine_cd_char} rules out affine contraction divisors on $\mathbb{R}$; nonlinear proper contraction divisors, however, do exist on $\mathbb{R}$.

\begin{example}\label{ex:nonlinear_1d}
Define continuous piecewise-linear mappings $f, g: \mathbb{R} \to \mathbb{R}$ by
\[
f(x) = \begin{cases} 2x, & \text{if } x \ge 0, \\ -\frac{x}{4}, & \text{if } x < 0, \end{cases}
\qquad
g(x) = \begin{cases} -\frac{x}{4}, & \text{if } x \ge 0, \\ 2x, & \text{if } x < 0. \end{cases}
\]
For every $n \in \mathbb{N}$, $|f^n(1) - f^n(2)| = 2^n \ge 2$ and $|g^n(-1) - g^n(-2)| = 2^n \ge 2$, so $f^n, g^n \notin \mathcal{C}(\mathbb{R})$ for all $n \in \mathbb{N}$.

Observe that $g(\mathbb{R}) = (-\infty, 0]$ and $f(\mathbb{R}) = [0, \infty)$. Since $f(u) = -u/4$ on $(-\infty, 0]$ and $g(u) = -u/4$ on $[0, \infty)$, composing $f$ and $g$ yields
\[
(f \circ g)(x) = \begin{cases} \frac{x}{16}, & \text{if } x \ge 0, \\ -\frac{x}{2}, & \text{if } x < 0, \end{cases}
\qquad
(g \circ f)(x) = \begin{cases} -\frac{x}{2}, & \text{if } x \ge 0, \\ \frac{x}{16}, & \text{if } x < 0. \end{cases}
\]
For any $x < 0 < y$, since $-x = |x| > 0$, we have
\[
|(f \circ g)(x) - (f \circ g)(y)| = \left|\frac{|x|}{2} - \frac{y}{16}\right| \le \max\left(\frac{|x|}{2}, \, \frac{y}{16}\right) \le \frac{1}{2}(|x| + y) = \frac{1}{2}|x - y|.
\]
Checking the remaining cases ($x, y \ge 0$ and $x, y \le 0$) gives $\Lip(f \circ g) = 1/2 < 1$; since $(g \circ f)(x) = -(f \circ g)(-x)$, it follows that $\Lip(g \circ f) = 1/2 < 1$ as well, so $(f, g)$ is a proper contraction divisor pair on $\mathbb{R}$.
\end{example}

In contrast to Theorem~\ref{thm:affine_cd_char}, left and right contraction divisors are strictly distinct among continuous nonlinear mappings.

\begin{example}\label{ex:one_sided_continuous}
On $X = [0, \infty)$ with the Euclidean metric, define continuous piecewise-linear mappings $f_1, g_1, f_2, g_2: X \to X$ by $g_1(x) = \min(2x, 1)$, $f_2(x) = 2x + 1$, and
\[
f_1(x) = \begin{cases} \frac{1}{4}(1 - x), & \text{if } 0 \le x \le 1, \\ 2(x - 1), & \text{if } x > 1, \end{cases}
\qquad
g_2(x) = \begin{cases} 2(1 - x), & \text{if } 0 \le x \le 1, \\ \frac{1}{4}(x - 1), & \text{if } x > 1. \end{cases}
\]
All four maps have Lipschitz constant $2$, so $f_1, g_1, f_2, g_2 \notin \mathcal{C}(X)$.

Since $g_1(X) = [0, 1]$, we have $(f_1 \circ g_1)(x) = \frac{1}{4}(1 - 2x)$ for $0 \le x \le 1/2$ and $(f_1 \circ g_1)(x) = 0$ for $x > 1/2$. Thus $\Lip(f_1 \circ g_1) = 1/2 < 1$, so $f_1$ is a left contraction divisor. However, $f_1$ is not a right contraction divisor: if $h: X \to X$ satisfies $h \circ f_1 \in \mathcal{C}(X)$ with contraction constant $k < 1$, then for any $u, v \in X$, setting $x = \frac{u}{2} + 1$ and $y = \frac{v}{2} + 1$ in $[1, \infty)$ gives $f_1(x) = u$, $f_1(y) = v$, and
\[
|h(u) - h(v)| = |h(f_1(x)) - h(f_1(y))| \le k|x - y| = \frac{k}{2}|u - v|,
\]
so $\Lip(h) \le k/2 < 1$ and $h \in \mathcal{C}(X)$.

On the other hand, $(g_2 \circ f_2)(x) = \frac{1}{2}x$ for all $x \in X$ has $\Lip(g_2 \circ f_2) = 1/2 < 1$, so $f_2$ is a right contraction divisor. Since $|f_2(x) - f_2(y)| = 2|x - y|$ for all $x, y \in X$, Proposition~\ref{prop:obstruction}(1) implies that $f_2$ is not a left contraction divisor (and shows that the surjectivity and density hypotheses in Proposition~\ref{prop:obstruction}(2) cannot be omitted).
\end{example}

We now establish the fixed-point and convergence properties of contraction divisors on complete metric spaces.

\begin{proposition}\label{prop:fixed_point}
Let $(X, d)$ be a complete metric space, and let $f, g: X \to X$ be mappings (not assumed to be continuous).
\begin{enumerate}
    \item If $f \circ g \in \mathcal{C}(X)$, then both $f \circ g$ and $g \circ f$ possess unique fixed points $x_0 \in X$ and $y_0 \in X$, respectively. Moreover, $(x_0, y_0)$ is the unique pair in $X \times X$ satisfying
    \begin{equation}\label{eq:periodic_cycle}
    g(x_0) = y_0 \quad \text{and} \quad f(y_0) = x_0.
    \end{equation}
    In particular, $x_0 = y_0$ is a common fixed point of $f$ and $g$ if and only if $g(x_0) = x_0$ (equivalently, $f(y_0) = y_0$), which holds whenever $f$ and $g$ commute.
    \item If $(f, g)$ is a contraction divisor pair, then for any initial point $x \in X$, the alternating sequence defined by $x_1 = g(x)$, $x_{2n} = f(x_{2n-1})$, and $x_{2n+1} = g(x_{2n})$ for $n \ge 1$ satisfies
    \[
    \lim_{n \to \infty} x_{2n} = x_0 \quad \text{and} \quad \lim_{n \to \infty} x_{2n+1} = y_0,
    \]
    where $x_0$ and $y_0$ are as in~(1).
\end{enumerate}
\end{proposition}

\begin{proof}
(1) Since $f \circ g$ is a contraction on the complete metric space $(X, d)$, Banach's Fixed Point Theorem gives a unique fixed point $x_0 \in X$ such that $f(g(x_0)) = x_0$. Set $y_0 = g(x_0)$. Then $f(y_0) = f(g(x_0)) = x_0$, and consequently
\[
(g \circ f)(y_0) = g(f(y_0)) = g(x_0) = y_0.
\]
Thus $y_0$ is a fixed point of $g \circ f$. To prove uniqueness, suppose $z \in X$ satisfies $(g \circ f)(z) = z$. Applying $f$ to both sides yields
\[
(f \circ g)(f(z)) = f(g(f(z))) = f(z),
\]
so $f(z)$ is a fixed point of $f \circ g$. By the uniqueness of $x_0$, we must have $f(z) = x_0$. Applying $g$ to both sides gives $z = g(f(z)) = g(x_0) = y_0$. Thus $y_0$ is unique. Furthermore, if any pair $(u, v) \in X \times X$ satisfies $g(u) = v$ and $f(v) = u$, then $(f \circ g)(u) = u$ and $(g \circ f)(v) = v$, forcing $u = x_0$ and $v = y_0$.

If $x_0 = y_0$ is a common fixed point of $f$ and $g$, then $g(x_0) = x_0$ and $f(y_0) = y_0$ trivially. Conversely, if $g(x_0) = x_0$, then \eqref{eq:periodic_cycle} yields $y_0 = g(x_0) = x_0$ and hence $f(x_0) = f(y_0) = x_0$; symmetrically, if $f(y_0) = y_0$, then $x_0 = f(y_0) = y_0$ and $g(y_0) = g(x_0) = y_0$. When $f$ and $g$ commute, $(g \circ f)(x_0) = (f \circ g)(x_0) = x_0$, so uniqueness gives $x_0 = y_0$, and \eqref{eq:periodic_cycle} reduces to $f(x_0) = g(x_0) = x_0$.

(2) By induction on the alternating sequence, $x_{2n} = (f \circ g)^n(x)$ and $x_{2n+1} = (g \circ f)^n(x_1)$ for all $n \ge 1$. Because $f \circ g$ and $g \circ f$ are contractions on the complete metric space $(X, d)$ with unique fixed points $x_0$ and $y_0$, respectively, Banach's Fixed Point Theorem yields $\lim_{n \to \infty} x_{2n} = x_0$ and $\lim_{n \to \infty} x_{2n+1} = y_0$.
\end{proof}

When $f \circ g \in \mathcal{C}(X)$, the even subsequence $x_{2n} = (f \circ g)^n(x)$ in Proposition~\ref{prop:fixed_point}(2) always converges to $x_0$; if $g$ is continuous at $x_0$, then $x_{2n+1} = g(x_{2n}) \to g(x_0) = y_0$ even without assuming $g \circ f \in \mathcal{C}(X)$. The following example shows that without continuity of $g$ at $x_0$ (when $f \circ g \in \mathcal{C}(X)$ but $g \circ f \notin \mathcal{C}(X)$), the odd subsequence $\{x_{2n+1}\}$ need not converge to $y_0$.

\begin{example}\label{ex:one_sided}
Define $f, g: \mathbb{R} \to \mathbb{R}$ by
\[
g(x) = \begin{cases} 0, & \text{if } x = 0, \\ \frac{1}{x}, & \text{if } x \neq 0, \end{cases}
\qquad
f(x) = \begin{cases} 0, & \text{if } x = 0, \\ \frac{1}{4x}, & \text{if } x \neq 0. \end{cases}
\]
Clearly $f, g \notin \mathcal{C}(\mathbb{R})$ (neither is continuous at $0$), while $(f \circ g)(x) = x/4$ is a contraction on $\mathbb{R}$, so $f$ is a left contraction divisor, and both $f \circ g$ and $g \circ f$ have the unique fixed point $x_0 = y_0 = 0$. However, $(g \circ f)(x) = 4x$ is not a contraction. Starting from $x = 1$, the alternating sequence gives $x_{2n} = 4^{-n} \to 0$, whereas $x_{2n+1} = 4^n \to \infty$.
\end{example}

\section{Coupled Contraction-Divisor Iterated Function Systems}\label{sec:cd_ifs}

Let $(X, d)$ be a complete metric space, and let $\mathcal{K}(X)$ denote the collection of all non-empty compact subsets of $X$ equipped with the Hausdorff metric $h_d$. For a classical IFS $\{f, g\}$, the attractor satisfies $A = f(A) \cup g(A)$. When $(f, g)$ is a proper contraction divisor pair, however, no non-empty compact set need satisfy this equation, as the following example shows.

\begin{example}\label{ex:no_union_attractor}
Let $f, g: \mathbb{R}^2 \to \mathbb{R}^2$ be defined by
\[
f(x, y) = \left(2x, \, \frac{y}{4}\right), \qquad g(x, y) = \left(\frac{x}{4} + 1, \, 2y\right).
\]
Here $(f \circ g)(x, y) = (x/2 + 2, y/2)$ and $(g \circ f)(x, y) = (x/2 + 1, y/2)$ are contractions with ratio $1/2$, admitting the unique fixed points $x_0 = (4, 0)$ and $y_0 = (2, 0)$, respectively. Suppose there exists $A \in \mathcal{K}(\mathbb{R}^2)$ satisfying $A = f(A) \cup g(A)$. Then $(f \circ g)(A) \subseteq f(A) \subseteq A$ and $(g \circ f)(A) \subseteq g(A) \subseteq A$. Choosing any $z \in A$, the sequences $\{(f \circ g)^n(z)\}$ and $\{(g \circ f)^n(z)\}$ remain in $A$ and converge to $x_0$ and $y_0$, respectively; since $A$ is closed, $x_0, y_0 \in A$. Because $f^n(A) \subseteq A$ for all $n \in \mathbb{N}$, we obtain
\[
d(f^n(x_0), f^n(y_0)) = d((4 \cdot 2^n, 0), \, (2 \cdot 2^n, 0)) = 2^{n+1} \to \infty \quad \text{as } n \to \infty,
\]
contradicting the boundedness of $A$.
\end{example}

To avoid repeated iterates of a single map, we follow the alternating structure of Proposition~\ref{prop:fixed_point} and consider a coupled system on $\mathcal{K}(X) \times \mathcal{K}(X)$ that combines contraction divisors with classical contractions.

\begin{definition}\label{def:coupled_cd_ifs}
Let $(X, d)$ be a complete metric space. Let $\mathcal{F} = \{f_1, \dots, f_{N_1}\}$ and $\tilde{\mathcal{F}} = \{\tilde{f}_1, \dots, \tilde{f}_{N_2}\}$ be non-empty finite collections of continuous mappings on $X$ such that $(f_i, \tilde{f}_j)$ is a contraction divisor pair for every $i \in \{1, \dots, N_1\}$ and $j \in \{1, \dots, N_2\}$, and let $\mathcal{G} = \{g_1, \dots, g_p\}$ and $\tilde{\mathcal{G}} = \{\tilde{g}_1, \dots, \tilde{g}_q\}$ be finite (possibly empty) collections of contractions in $\mathcal{C}(X)$. The four-tuple $(\mathcal{F}, \tilde{\mathcal{F}}, \mathcal{G}, \tilde{\mathcal{G}})$ is called a \emph{Coupled Contraction-Divisor Iterated Function System (CD-IFS)}. The \emph{coupled Hutchinson operator} $\mathbf{T}: \mathcal{K}(X) \times \mathcal{K}(X) \to \mathcal{K}(X) \times \mathcal{K}(X)$ is defined by $\mathbf{T}(B_1, B_2) = (T_1(B_1, B_2), T_2(B_1, B_2))$, where
\[
T_1(B_1, B_2) = \bigcup_{i=1}^{N_1} f_i(B_2) \cup \bigcup_{k=1}^p g_k(B_1), \qquad T_2(B_1, B_2) = \bigcup_{j=1}^{N_2} \tilde{f}_j(B_1) \cup \bigcup_{l=1}^q \tilde{g}_l(B_2).
\]
A pair $(A_1, A_2) \in \mathcal{K}(X) \times \mathcal{K}(X)$ is called a \emph{coupled attractor} of the system if $\mathbf{T}(A_1, A_2) = (A_1, A_2)$, that is,
\begin{equation}\label{eq:coupled_invariance}
A_1 = \bigcup_{i=1}^{N_1} f_i(A_2) \cup \bigcup_{k=1}^p g_k(A_1), \qquad A_2 = \bigcup_{j=1}^{N_2} \tilde{f}_j(A_1) \cup \bigcup_{l=1}^q \tilde{g}_l(A_2).
\end{equation}
\end{definition}

When at least one of the contraction families $\mathcal{G}$ or $\tilde{\mathcal{G}}$ is empty, \eqref{eq:coupled_invariance} decouples into a finite classical IFS.

\begin{proposition}\label{prop:coupled_existence}
Let $(\mathcal{F}, \tilde{\mathcal{F}}, \mathcal{G}, \tilde{\mathcal{G}})$ be a CD-IFS on a complete metric space $(X, d)$. If $\tilde{\mathcal{G}} = \emptyset$, then the system admits a unique coupled attractor $(A_1, A_2) \in \mathcal{K}(X) \times \mathcal{K}(X)$, where $A_1$ is the classical attractor of the finite contractive IFS $\{f_i \circ \tilde{f}_j : 1 \le i \le N_1, \, 1 \le j \le N_2\} \cup \mathcal{G}$ and $A_2 = \bigcup_{j=1}^{N_2} \tilde{f}_j(A_1)$. Symmetrically, if $\mathcal{G} = \emptyset$, there is a unique coupled attractor $(A_1, A_2)$, where $A_2$ is the classical attractor of $\{\tilde{f}_j \circ f_i : 1 \le i \le N_1, \, 1 \le j \le N_2\} \cup \tilde{\mathcal{G}}$ and $A_1 = \bigcup_{i=1}^{N_1} f_i(A_2)$.
\end{proposition}

\begin{proof}
Suppose $\tilde{\mathcal{G}} = \emptyset$; the case $\mathcal{G} = \emptyset$ is symmetric. Substituting $A_2 = \bigcup_{j=1}^{N_2} \tilde{f}_j(A_1)$ into the first equation of \eqref{eq:coupled_invariance} shows that $(A_1, A_2) \in \mathcal{K}(X) \times \mathcal{K}(X)$ satisfies \eqref{eq:coupled_invariance} if and only if
\[
A_1 = \bigcup_{i=1}^{N_1} \bigcup_{j=1}^{N_2} (f_i \circ \tilde{f}_j)(A_1) \cup \bigcup_{k=1}^p g_k(A_1) \quad \text{and} \quad A_2 = \bigcup_{j=1}^{N_2} \tilde{f}_j(A_1).
\]
Since every $f_i \circ \tilde{f}_j$ and every $g_k$ belongs to $\mathcal{C}(X)$, Hutchinson's theorem \cite{Hutchinson1981} gives a unique $A_1 \in \mathcal{K}(X)$ satisfying the first equality. The second equality then uniquely determines $A_2 = \bigcup_{j=1}^{N_2} \tilde{f}_j(A_1)$, which belongs to $\mathcal{K}(X)$ because each $\tilde{f}_j$ is continuous and $\tilde{\mathcal{F}}$ is non-empty.
\end{proof}

When both $\mathcal{G}$ and $\tilde{\mathcal{G}}$ are non-empty, \eqref{eq:coupled_invariance} no longer decouples into a finite classical IFS. Moreover, without a compatibility condition between the contraction divisors and the classical contractions, a coupled attractor need not exist, even on $\mathbb{R}^2$ when all compositions $f_i \circ \tilde{f}_j$ and $\tilde{f}_j \circ f_i$ as well as all maps in $\mathcal{G} \cup \tilde{\mathcal{G}}$ are similarity contractions.

\begin{example}\label{ex:no_coupled_attractor}
Consider the CD-IFS $(\{f_1\}, \{\tilde{f}_1\}, \{g_1\}, \{\tilde{g}_1\})$ on $\mathbb{R}^2$ defined by
\begin{gather*}
f_1(x, y) = \left(4x, \, \frac{y}{8}\right), \qquad \tilde{f}_1(x, y) = \left(\frac{x}{8}, \, 4y\right), \\
g_1(x, y) = \left(\frac{y}{2} + 1, \, \frac{x}{2}\right), \qquad \tilde{g}_1(x, y) = \left(\frac{y}{2}, \, \frac{x}{2}\right).
\end{gather*}
Here $\Lip(f_1^n) = \Lip(\tilde{f}_1^n) = 4^n > 1$ for all $n \in \mathbb{N}$, while $(f_1 \circ \tilde{f}_1)(x, y) = (\tilde{f}_1 \circ f_1)(x, y) = (x/2, y/2)$ is a similarity contraction with ratio $1/2$; thus $(f_1, \tilde{f}_1)$ is a proper contraction divisor pair. Moreover, $g_1$ and $\tilde{g}_1$ are similarity contractions with ratio $1/2$.

Suppose there exists a coupled attractor $(A_1, A_2) \in \mathcal{K}(\mathbb{R}^2) \times \mathcal{K}(\mathbb{R}^2)$ satisfying $A_1 = f_1(A_2) \cup g_1(A_1)$ and $A_2 = \tilde{f}_1(A_1) \cup \tilde{g}_1(A_2)$. Then $(f_1 \circ \tilde{f}_1)(A_1) \subseteq f_1(A_2) \subseteq A_1$. Choosing any $z \in A_1$, the sequence $\{(f_1 \circ \tilde{f}_1)^n(z)\} = \{2^{-n} z\}$ remains in $A_1$ and converges to $(0, 0)$; since $A_1$ is closed, $(0, 0) \in A_1$. Because $g_1(A_1) \subseteq A_1$, it follows that $(1, 0) = g_1(0, 0) \in A_1$. Consider the composition $\Psi = f_1 \circ \tilde{g}_1 \circ \tilde{f}_1 \circ g_1$, which satisfies $\Psi(A_1) \subseteq A_1$ and is given by
\[
\Psi(x, y) = \left(4x, \, \frac{y}{256} + \frac{1}{128}\right).
\]
For every $n \in \mathbb{N}$, the point $\Psi^n(1, 0) \in A_1$ has first coordinate $4^n$, so $\|\Psi^n(1, 0)\| \ge 4^n$; letting $n \to \infty$ contradicts the boundedness of $A_1$. Hence no coupled attractor exists.
\end{example}

In Example~\ref{ex:no_coupled_attractor}, the linear parts of $g_1$ and $\tilde{g}_1$ interchange the expanding and contracting directions of $f_1$ and $\tilde{f}_1$, producing an expanding composition. In Section~\ref{sec:dim}, we impose an alignment condition on $\mathbb{R}^m$ that rules out this behavior and yields both a unique coupled attractor and an exact Hausdorff dimension formula.

\section{Hausdorff Dimension of Coupled Self-Affine Sets}\label{sec:dim}

In this section, we generalize the dimension theory of \cite{AmalSelfAffine2026} to coupled self-affine systems on $\mathbb{R}^m$ combining affine contraction divisors with similarity contractions. Throughout, we equip $\mathcal{K}(\mathbb{R}^m) \times \mathcal{K}(\mathbb{R}^m)$ with the product Hausdorff metric
\[
H((B_1, B_2), (C_1, C_2)) = \max\{h_d(B_1, C_1), \, h_d(B_2, C_2)\}.
\]
We first motivate the general dimension formula with an explicit planar example in which the attractor decomposes into countably many self-similar copies of itself.

\begin{figure}[htbp]
    \centering
    \includegraphics[width=\textwidth]{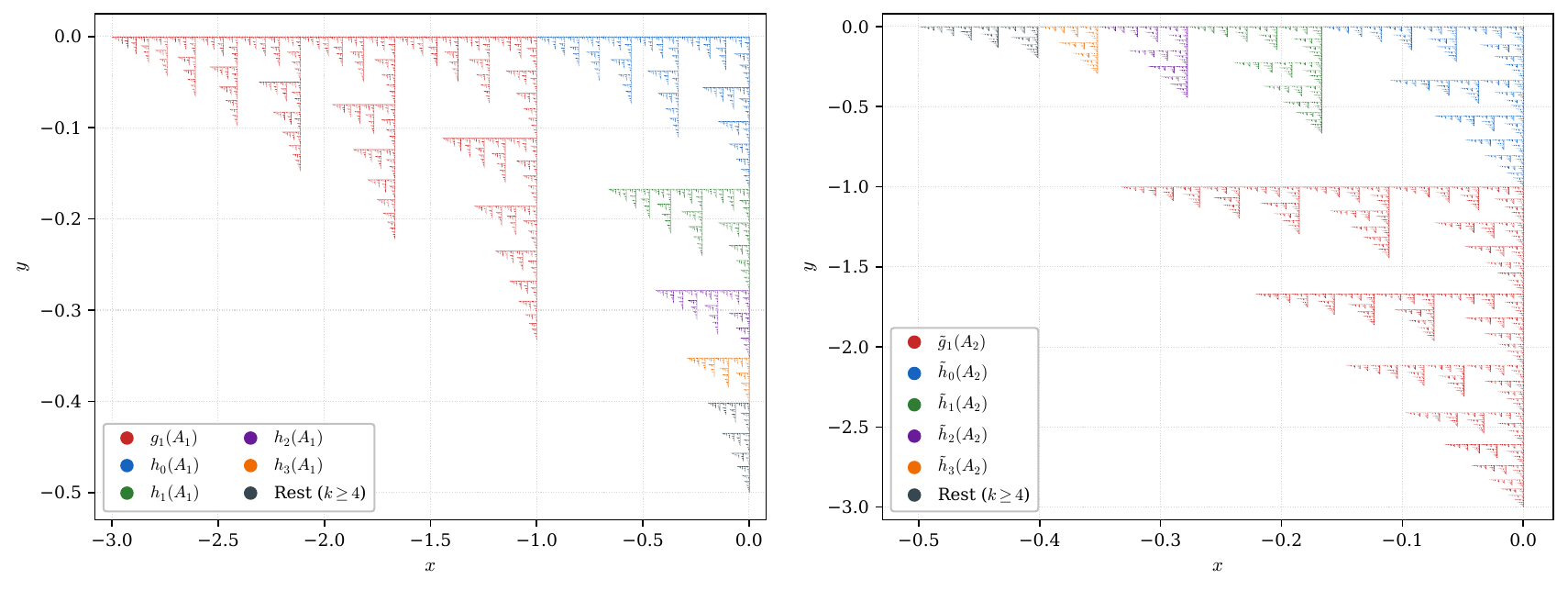}
    \caption{Coupled attractors $A_1$ (left) and $A_2$ (right) of Example~\ref{ex:motivating_sandwich} decomposed into their countable self-similar components $g_1(A_1)$, $h_k(A_1) = (f_1 \circ \tilde{g}_1^k \circ \tilde{f}_1)(A_1)$ and $\tilde{g}_1(A_2)$, $\tilde{h}_k(A_2) = (\tilde{f}_1 \circ g_1^k \circ f_1)(A_2)$ ($k = 0, 1, 2, 3$ and $k \ge 4$).}
    \label{fig:sandwich}
\end{figure}

\begin{example}\label{ex:motivating_sandwich}
Consider the CD-IFS $(\{f_1\}, \{\tilde{f}_1\}, \{g_1\}, \{\tilde{g}_1\})$ on $\mathbb{R}^2$ given by
\begin{gather*}
f_1(x, y) = \left(2x, \, \frac{y}{6}\right), \qquad \tilde{f}_1(x, y) = \left(\frac{x}{6}, \, 2y\right), \\
g_1(x, y) = \left(\frac{2}{3}x - 1, \, \frac{2}{3}y\right), \qquad \tilde{g}_1(x, y) = \left(\frac{2}{3}x, \, \frac{2}{3}y - 1\right).
\end{gather*}
Neither $f_1$ nor $\tilde{f}_1$ is a contraction, since $f_1$ expands the first coordinate by a factor of $2$ and $\tilde{f}_1$ expands the second by $2$. However, $(f_1 \circ \tilde{f}_1)(x, y) = (\tilde{f}_1 \circ f_1)(x, y) = (x/3, y/3)$ is a similarity with contraction ratio $1/3$, and $g_1, \tilde{g}_1$ are similarities with ratio $2/3$. Expanding $\mathbf{T}^3(B_1, B_2)$ yields eight three-fold compositions in each coordinate; since the linear parts of $g_1$ and $\tilde{g}_1$ are $\frac{2}{3}I_2$, each of these compositions has Lipschitz constant at most $8/9 < 1$. Thus $\mathbf{T}^3$ is a contraction on $(\mathcal{K}(\mathbb{R}^2) \times \mathcal{K}(\mathbb{R}^2), H)$, so $\mathbf{T}^3$ has a unique fixed point $(A_1, A_2)$, which is therefore the unique coupled attractor of the system.

Consider the open rectangles $U_1 = (-3, 0) \times (-1/2, 0)$ and $U_2 = (-1/2, 0) \times (-3, 0)$. We compute
\[
f_1(U_2) = (-1, 0) \times \left(-\frac{1}{2}, 0\right) \subset U_1, \qquad g_1(U_1) = (-3, -1) \times \left(-\frac{1}{3}, 0\right) \subset U_1.
\]
Since their $x$-projections $(-1, 0)$ and $(-3, -1)$ are disjoint, $f_1(U_2) \cap g_1(U_1) = \emptyset$. Symmetrically, $\tilde{f}_1(U_1) = (-1/2, 0) \times (-1, 0)$ and $\tilde{g}_1(U_2) = (-1/3, 0) \times (-3, -1)$ are disjoint subsets of $U_2$.

Iterating $A_2 = \tilde{g}_1(A_2) \cup \tilde{f}_1(A_1)$ $n$ times yields $A_2 = \tilde{g}_1^n(A_2) \cup \bigcup_{k=0}^{n-1} (\tilde{g}_1^k \circ \tilde{f}_1)(A_1)$. Because $\tilde{g}_1(A_2) \subseteq A_2$ and $A_2$ is non-empty and compact, the unique fixed point $(0, -3)$ of $\tilde{g}_1$ belongs to $A_2$, and $\bigcap_{n=1}^\infty \tilde{g}_1^n(A_2) = \{(0, -3)\}$. Any point of $A_2$ outside $\bigcup_{k=0}^\infty (\tilde{g}_1^k \circ \tilde{f}_1)(A_1)$ must lie in $\tilde{g}_1^n(A_2)$ for every $n \ge 1$, so $A_2 = \bigcup_{k=0}^\infty (\tilde{g}_1^k \circ \tilde{f}_1)(A_1) \cup \{(0, -3)\}$. Substituting this into $A_1 = g_1(A_1) \cup f_1(A_2)$ gives the countable decomposition
\begin{equation}\label{eq:motivating_decomp}
A_1 = g_1(A_1) \cup \bigcup_{k=0}^\infty h_k(A_1) \cup \left\{\left(0, -\frac{1}{2}\right)\right\},
\end{equation}
where $h_k = f_1 \circ \tilde{g}_1^k \circ \tilde{f}_1$ and the singleton $f_1(0, -3) = (0, -1/2)$ has $\mathcal{H}^s$-measure zero for $s > 0$. Since $\tilde{g}_1^k(u, v) = ((2/3)^k u, (2/3)^k v - c_k)$ with $c_k = \sum_{j=0}^{k-1} (2/3)^j = 3(1 - (2/3)^k)$, composing the maps gives
\[
h_k(x, y) = \left(\frac{1}{3}\left(\frac{2}{3}\right)^k x, \, \frac{1}{3}\left(\frac{2}{3}\right)^k y - \frac{c_k}{6}\right),
\]
which is a similarity with contraction ratio $r_k = \frac{1}{3}(2/3)^k$. Because $\tilde{f}_1(U_1)$ and $\tilde{g}_1(U_2)$ are disjoint subsets of $U_2$, applying $\tilde{g}_1$ inductively shows that the sets $\{\tilde{g}_1^k(\tilde{f}_1(U_1))\}_{k=0}^\infty$ are pairwise disjoint inside $U_2$; applying the injective map $f_1$ then yields pairwise disjoint open sets $h_k(U_1) \subset f_1(U_2)$ ($k \ge 0$), all of which are disjoint from $g_1(U_1)$. Proceeding heuristically under the assumptions that the components in \eqref{eq:motivating_decomp} are $\mathcal{H}^s$-almost disjoint and $0 < \mathcal{H}^s(A_1) < \infty$, we apply $\mathcal{H}^s$ to \eqref{eq:motivating_decomp} and divide by $\mathcal{H}^s(A_1)$ to obtain
\[
1 = \left(\frac{2}{3}\right)^s + \sum_{k=0}^\infty \left(\frac{1}{3}\left(\frac{2}{3}\right)^k\right)^s = \left(\frac{2}{3}\right)^s + \frac{(1/3)^s}{1 - (2/3)^s}.
\]
Subtracting $(2/3)^s$ from both sides and multiplying by $1 - (2/3)^s$ gives $\bigl(1 - (2/3)^s\bigr)^2 = (1/3)^s$, and taking the positive square root yields $(2/3)^s + (1/\sqrt{3})^s = 1$, whence $s \approx 1.4635$. The coupled attractors $(A_1, A_2)$ and their countable self-similar decompositions are illustrated in Figure~\ref{fig:sandwich}.
\end{example}

The essential feature of Example~\ref{ex:motivating_sandwich} is that every composition $f_1 \circ \tilde{g}_1^k \circ \tilde{f}_1$ reduces to a similarity contraction. We now formulate the general conditions on $\mathbb{R}^m$ that preserve this property and yield a rigorous existence and dimension theory.

\begin{definition}\label{def:sim_cd}
Let $f(x) = M x + a$ and $\tilde{f}(x) = \tilde{M} x + b$ be non-contractive affine maps on $\mathbb{R}^m$. The pair $(f, \tilde{f})$ is called a \emph{similarity contraction divisor pair} if both $f \circ \tilde{f}$ and $\tilde{f} \circ f$ are similarity contractions on $\mathbb{R}^m$; that is, there exists $c \in (0, 1)$ such that
\[
\|(f \circ \tilde{f})(x) - (f \circ \tilde{f})(y)\| = c\|x - y\| \quad \text{and} \quad \|(\tilde{f} \circ f)(x) - (\tilde{f} \circ f)(y)\| = c\|x - y\|
\]
for all $x, y \in \mathbb{R}^m$.
\end{definition}

\begin{definition}[{\cite[Definition 1.2 and Definition 5.1]{AmalSelfAffine2026}}]\label{def:aligned_sim}
Let $f: \mathbb{R}^m \to \mathbb{R}^m$ be an affine transformation with linear part $M$. A similarity transformation $g: \mathbb{R}^m \to \mathbb{R}^m$ with linear part $S$ is called \emph{$f$-aligned} if
\[
S(M^\top M) = (M^\top M)S.
\]
Given a finite collection of affine maps $\mathcal{F} = \{f_1, \dots, f_{N_1}\}$ on $\mathbb{R}^m$, the similarity $g$ is said to be \emph{universally aligned} with $\mathcal{F}$ if $g$ is $f_i$-aligned for every $i \in \{1, \dots, N_1\}$.
\end{definition}

Consider a CD-IFS $(\mathcal{F}, \tilde{\mathcal{F}}, \mathcal{G}, \tilde{\mathcal{G}})$ on $\mathbb{R}^m$, where $\mathcal{F} = \{f_1, \dots, f_{N_1}\}$ and $\tilde{\mathcal{F}} = \{\tilde{f}_1, \dots, \tilde{f}_{N_2}\}$ consist of invertible affine maps, and $\mathcal{G} = \{g_1, \dots, g_p\}$ and $\tilde{\mathcal{G}} = \{\tilde{g}_1, \dots, \tilde{g}_q\}$ consist of similarity contractions with contraction ratios $r_k \in (0, 1)$ ($1 \le k \le p$) and $\tilde{r}_l \in (0, 1)$ ($1 \le l \le q$), respectively. To describe the $L$-th iterate $\mathbf{T}^L(B_1, B_2) = (T_1^{(L)}(B_1, B_2), T_2^{(L)}(B_1, B_2))$ for $L \in \mathbb{N}$, we refer to the first and second coordinates as State~$1$ and State~$2$, respectively. We treat the four collections $\mathcal{F}$, $\tilde{\mathcal{F}}$, $\mathcal{G}$, and $\tilde{\mathcal{G}}$ as formally disjoint and call the elements of $\mathcal{F} \cup \tilde{\mathcal{F}} \cup \mathcal{G} \cup \tilde{\mathcal{G}}$ \emph{symbols}; each symbol $\psi$ thus uniquely identifies the collection to which it belongs, even if maps from different collections coincide as functions on $\mathbb{R}^m$. In the coupled Hutchinson operator $\mathbf{T}$, the first coordinate $T_1$ applies symbols from $\mathcal{F} \cup \mathcal{G}$, whereas the second coordinate $T_2$ applies symbols from $\tilde{\mathcal{F}} \cup \tilde{\mathcal{G}}$. To record which coordinate each symbol acts on, define $\tau(\psi) \in \{1, 2\}$ by
\[
\tau(\psi) = \begin{cases} 1, & \text{if } \psi \in \mathcal{G} \cup \tilde{\mathcal{F}}, \\ 2, & \text{if } \psi \in \mathcal{F} \cup \tilde{\mathcal{G}}, \end{cases}
\]
so that
\[
T_1(B_1, B_2) = \bigcup_{\psi \in \mathcal{F} \cup \mathcal{G}} \psi(B_{\tau(\psi)}), \qquad T_2(B_1, B_2) = \bigcup_{\psi \in \tilde{\mathcal{F}} \cup \tilde{\mathcal{G}}} \psi(B_{\tau(\psi)}).
\]
Accordingly, a finite sequence $\omega = (\psi_1, \dots, \psi_L)$ of symbols from $\mathcal{F} \cup \tilde{\mathcal{F}} \cup \mathcal{G} \cup \tilde{\mathcal{G}}$ is called an \emph{admissible word of length $L$ starting from State~$1$} if:
\begin{enumerate}
    \item $\psi_1 \in \mathcal{F} \cup \mathcal{G}$;
    \item for each $1 \le t < L$, $\psi_{t+1} \in \mathcal{F} \cup \mathcal{G}$ if $\tau(\psi_t) = 1$, and $\psi_{t+1} \in \tilde{\mathcal{F}} \cup \tilde{\mathcal{G}}$ if $\tau(\psi_t) = 2$.
\end{enumerate}
An \emph{admissible word of length $L$ starting from State~$2$} is defined identically, except that condition~(1) is replaced by $\psi_1 \in \tilde{\mathcal{F}} \cup \tilde{\mathcal{G}}$. We denote by $\mathcal{I}_1^L$ and $\mathcal{I}_2^L$ the sets of all admissible words of length $L$ starting from State~$1$ and State~$2$, respectively. For each $\omega = (\psi_1, \dots, \psi_L) \in \mathcal{I}_1^L \cup \mathcal{I}_2^L$, set $h_\omega = \psi_1 \circ \dots \circ \psi_L$ and extend $\tau$ to admissible words by setting $\tau(\omega) = \tau(\psi_L)$, the \emph{terminal state index} of $\omega$. By induction on $L$, substituting $T_1(B_1, B_2)$ for $B_1$ and $T_2(B_1, B_2)$ for $B_2$ at each step yields
\[
T_1^{(L)}(B_1, B_2) = \bigcup_{\omega \in \mathcal{I}_1^L} h_\omega(B_{\tau(\omega)}), \qquad T_2^{(L)}(B_1, B_2) = \bigcup_{\omega \in \mathcal{I}_2^L} h_\omega(B_{\tau(\omega)}).
\]
In particular, by condition~(2), symbols from $\mathcal{F}$ and $\tilde{\mathcal{F}}$ alternate in every admissible word $\omega \in \mathcal{I}_1^L \cup \mathcal{I}_2^L$: after any $f_i \in \mathcal{F}$, only symbols from $\tilde{\mathcal{G}}$ may appear before either a symbol from $\tilde{\mathcal{F}}$ occurs or the word terminates; symmetrically, after any $\tilde{f}_j \in \tilde{\mathcal{F}}$, only symbols from $\mathcal{G}$ may appear before either a symbol from $\mathcal{F}$ occurs or the word terminates.

\begin{lemma}\label{lem:sandwich}
Let $\mathcal{F} = \{f_1, \dots, f_{N_1}\}$ and $\tilde{\mathcal{F}} = \{\tilde{f}_1, \dots, \tilde{f}_{N_2}\}$ be non-empty finite sets of invertible affine maps on $\mathbb{R}^m$ with linear parts $M_i$ ($1 \le i \le N_1$) and $\tilde{M}_j$ ($1 \le j \le N_2$), respectively, and let $\mathcal{G} = \{g_1, \dots, g_p\}$ and $\tilde{\mathcal{G}} = \{\tilde{g}_1, \dots, \tilde{g}_q\}$ be finite (possibly empty) sets of similarity contractions with contraction ratios $r_k \in (0, 1)$ ($1 \le k \le p$) and $\tilde{r}_l \in (0, 1)$ ($1 \le l \le q$), respectively. Assume that every pair $(f_i, \tilde{f}_j)$ is a similarity contraction divisor pair whose compositions $f_i \circ \tilde{f}_j$ and $\tilde{f}_j \circ f_i$ have contraction ratio $c_{i,j} \in (0, 1)$, every map in $\mathcal{G}$ is universally aligned with $\tilde{\mathcal{F}}$, and every map in $\tilde{\mathcal{G}}$ is universally aligned with $\mathcal{F}$. Then:
\begin{enumerate}
    \item For any $1 \le i \le N_1$, $1 \le j \le N_2$, any integer $d \ge 0$, and any maps $\tilde{g}_{l_1}, \dots, \tilde{g}_{l_d} \in \tilde{\mathcal{G}}$, the composition
    \[
    \Phi = f_i \circ \tilde{g}_{l_1} \circ \dots \circ \tilde{g}_{l_d} \circ \tilde{f}_j
    \]
    is a similarity contraction on $\mathbb{R}^m$ with contraction ratio $c_{i,j} \tilde{r}_{l_1} \cdots \tilde{r}_{l_d}$. Symmetrically, for any $g_{k_1}, \dots, g_{k_d} \in \mathcal{G}$, the composition $\tilde{\Phi} = \tilde{f}_j \circ g_{k_1} \circ \dots \circ g_{k_d} \circ f_i$ is a similarity contraction with ratio $c_{i,j} r_{k_1} \cdots r_{k_d}$.
    \item Let $\mathrm{Id}$ denote the identity map on $\mathbb{R}^m$. For any $L \in \mathbb{N}$, every admissible composition $h_\omega$ ($\omega \in \mathcal{I}_1^L$) factors as
    \begin{equation}\label{eq:inplace_factor}
    h_\omega = h^{\mathrm{head}} \circ h^{\mathrm{tail}} \circ h^{\mathrm{sim}},
    \end{equation}
    where $h^{\mathrm{head}}$ is a composition of zero or more maps from $\mathcal{G}$ and blocks of the form $\Phi$, $h^{\mathrm{tail}} \in \mathcal{T} := \{\mathrm{Id}\} \cup \mathcal{F}$, and $h^{\mathrm{sim}}$ is a composition of zero or more maps from $\tilde{\mathcal{G}}$ (with an empty composition understood as $\mathrm{Id}$). In particular, both $h^{\mathrm{head}}$ and $h^{\mathrm{sim}}$ are similarities; if $\tau(\omega) = 1$, then $h^{\mathrm{tail}} = h^{\mathrm{sim}} = \mathrm{Id}$, while if $\tau(\omega) = 2$, then $h^{\mathrm{tail}} \in \mathcal{F}$. A symmetric factorization holds for $\omega \in \mathcal{I}_2^L$ with $h^{\mathrm{tail}} \in \{\mathrm{Id}\} \cup \tilde{\mathcal{F}}$.
\end{enumerate}
\end{lemma}

\begin{proof}
(1) Let $\tilde{S}_l$ denote the linear part of $\tilde{g}_l$, so $\tilde{S}_l^\top \tilde{S}_l = \tilde{r}_l^2 I_m$. Since $\tilde{g}_l$ is $f_i$-aligned, $\tilde{S}_l(M_i^\top M_i) = (M_i^\top M_i)\tilde{S}_l$, which yields
\[
\tilde{S}_l^\top (M_i^\top M_i) \tilde{S}_l = \tilde{S}_l^\top \tilde{S}_l (M_i^\top M_i) = \tilde{r}_l^2 M_i^\top M_i.
\]
Setting $\tilde{S} = \tilde{S}_{l_1} \cdots \tilde{S}_{l_d}$ (with $\tilde{S} = I_m$ when $d = 0$), induction on $d$ gives $\tilde{S}^\top (M_i^\top M_i) \tilde{S} = (\tilde{r}_{l_1} \cdots \tilde{r}_{l_d})^2 M_i^\top M_i$. Since $f_i \circ \tilde{f}_j$ is a similarity with ratio $c_{i,j}$, its linear part $M_i \tilde{M}_j$ satisfies $(M_i \tilde{M}_j)^\top (M_i \tilde{M}_j) = \tilde{M}_j^\top (M_i^\top M_i) \tilde{M}_j = c_{i,j}^2 I_m$. Therefore, the linear part $Q = M_i \tilde{S} \tilde{M}_j$ of $\Phi$ satisfies
\[
Q^\top Q = \tilde{M}_j^\top \tilde{S}^\top (M_i^\top M_i) \tilde{S} \tilde{M}_j = (\tilde{r}_{l_1} \cdots \tilde{r}_{l_d})^2 \tilde{M}_j^\top (M_i^\top M_i) \tilde{M}_j = (c_{i,j} \tilde{r}_{l_1} \cdots \tilde{r}_{l_d})^2 I_m.
\]
Hence $\|Q v\| = (c_{i,j} \tilde{r}_{l_1} \cdots \tilde{r}_{l_d})\|v\|$ for all $v \in \mathbb{R}^m$, proving that $\Phi$ is a similarity contraction with ratio $c_{i,j} \tilde{r}_{l_1} \cdots \tilde{r}_{l_d}$. The proof for $\tilde{\Phi}$ is symmetric.

(2) By the transition rules for $\omega \in \mathcal{I}_1^L$, the maps from $\mathcal{F}$ and $\tilde{\mathcal{F}}$ alternate. If $\tau(\omega) = 1$, then $\omega$ contains an equal number of maps from $\mathcal{F}$ and $\tilde{\mathcal{F}}$, so $h_\omega$ is a composition of maps from $\mathcal{G}$ and blocks of the form $\Phi = f_i \circ \tilde{g}_{l_1} \circ \dots \circ \tilde{g}_{l_d} \circ \tilde{f}_j$; setting $h^{\mathrm{head}} = h_\omega$ and $h^{\mathrm{tail}} = h^{\mathrm{sim}} = \mathrm{Id}$ gives \eqref{eq:inplace_factor}. If $\tau(\omega) = 2$, then $\omega$ ends with a final unmatched map $f_{i^*} \in \mathcal{F}$ followed by zero or more maps from $\tilde{\mathcal{G}}$; letting $h^{\mathrm{head}}$ be the composition of all maps preceding $f_{i^*}$ (with $h^{\mathrm{head}} = \mathrm{Id}$ if $f_{i^*}$ is the first map in $\omega$), setting $h^{\mathrm{tail}} = f_{i^*} \in \mathcal{F}$, and setting $h^{\mathrm{sim}}$ to the trailing composition from $\tilde{\mathcal{G}}$ yields \eqref{eq:inplace_factor}. By part (1), $h^{\mathrm{head}}$ is a composition of similarities, hence a similarity.
\end{proof}

\begin{theorem}\label{thm:hybrid_coupled_existence}
Under the hypotheses of Lemma~\ref{lem:sandwich}, the CD-IFS $(\mathcal{F}, \tilde{\mathcal{F}}, \mathcal{G}, \tilde{\mathcal{G}})$ admits a unique coupled attractor $(A_1, A_2) \in \mathcal{K}(\mathbb{R}^m) \times \mathcal{K}(\mathbb{R}^m)$. Moreover, for every initial pair $(B_1, B_2) \in \mathcal{K}(\mathbb{R}^m) \times \mathcal{K}(\mathbb{R}^m)$, the iterates $\mathbf{T}^n(B_1, B_2)$ converge to $(A_1, A_2)$ in the product Hausdorff metric $H$ as $n \to \infty$.
\end{theorem}

\begin{proof}
Set $c_{\max} = \max_{i,j} c_{i,j} \in (0, 1)$, $\theta = \sqrt{c_{\max}} \in (0, 1)$, and
\[
M_{\max} = \max\big(\{1\} \cup \{\|M_i\| : 1 \le i \le N_1\} \cup \{\|\tilde{M}_j\| : 1 \le j \le N_2\}\big) < \infty.
\]
Set $r_{\max} = 0$ if $p = q = 0$ (with the convention $0^0 = 1$) and $r_{\max} = \max(\{r_k : 1 \le k \le p\} \cup \{\tilde{r}_l : 1 \le l \le q\}) \in (0, 1)$ if $p + q \ge 1$, and let $\rho = \max\{\theta, r_{\max}\} \in (0, 1)$.

Fix $L \in \mathbb{N}$ and $\omega \in \mathcal{I}_1^L$. Suppose $\omega$ contains $m_{\mathrm{sim}}$ symbols from $\mathcal{G} \cup \tilde{\mathcal{G}}$ and $k = L - m_{\mathrm{sim}}$ symbols from $\mathcal{F} \cup \tilde{\mathcal{F}}$. In the factorization $h_\omega = h^{\mathrm{head}} \circ h^{\mathrm{tail}} \circ h^{\mathrm{sim}}$ of Lemma~\ref{lem:sandwich}(2), $h^{\mathrm{head}}$ contains $\lfloor k/2 \rfloor$ blocks of the form $\Phi$, each having similarity ratio at most $c_{\max} = \theta^2 \le \rho^2$, every symbol from $\mathcal{G} \cup \tilde{\mathcal{G}}$ contributes a factor at most $r_{\max} \le \rho$, and $\Lip(h^{\mathrm{tail}}) \le M_{\max}$. Since $m_{\mathrm{sim}} + 2\lfloor k/2 \rfloor \ge L - 1$, it follows that
\[
\Lip(h_\omega) \le \Lip(h^{\mathrm{head}}) \Lip(h^{\mathrm{tail}}) \Lip(h^{\mathrm{sim}}) \le M_{\max} (r_{\max})^{m_{\mathrm{sim}}} (c_{\max})^{\lfloor k/2 \rfloor} \le M_{\max} \rho^{L - 1}.
\]
By symmetry, the same bound holds for every $\omega \in \mathcal{I}_2^L$. Consequently, for all $(B_1, B_2)$ and $(C_1, C_2)$ in $\mathcal{K}(\mathbb{R}^m) \times \mathcal{K}(\mathbb{R}^m)$ and every $L \in \mathbb{N}$,
\begin{align*}
h_d(T_1^{(L)}(B_1, B_2), T_1^{(L)}(C_1, C_2)) &\le \max_{\omega \in \mathcal{I}_1^L} h_d(h_\omega(B_{\tau(\omega)}), h_\omega(C_{\tau(\omega)})) \\
&\le M_{\max} \rho^{L - 1} H((B_1, B_2), (C_1, C_2)),
\end{align*}
and symmetrically for $T_2^{(L)}$, so
\begin{equation}\label{eq:TL_bound}
H(\mathbf{T}^L(B_1, B_2), \mathbf{T}^L(C_1, C_2)) \le M_{\max} \rho^{L - 1} H((B_1, B_2), (C_1, C_2)).
\end{equation}

Since $\rho \in (0, 1)$, we may choose $L_0 \in \mathbb{N}$ such that $M_{\max} \rho^{L_0 - 1} < 1$. Because $(\mathcal{K}(\mathbb{R}^m) \times \mathcal{K}(\mathbb{R}^m), H)$ is complete, Banach's Fixed Point Theorem yields a unique fixed point $(A_1, A_2)$ of $\mathbf{T}^{L_0}$. Since $\mathbf{T}^{L_0}(\mathbf{T}(A_1, A_2)) = \mathbf{T}(\mathbf{T}^{L_0}(A_1, A_2)) = \mathbf{T}(A_1, A_2)$, the uniqueness of the fixed point implies $\mathbf{T}(A_1, A_2) = (A_1, A_2)$, and since every fixed point of $\mathbf{T}$ is also a fixed point of $\mathbf{T}^{L_0}$, $(A_1, A_2)$ is the unique coupled attractor. Finally, for any initial pair $(B_1, B_2)$, applying \eqref{eq:TL_bound} with $(C_1, C_2) = (A_1, A_2)$ gives
\begin{align*}
H(\mathbf{T}^n(B_1, B_2), (A_1, A_2)) &= H(\mathbf{T}^n(B_1, B_2), \mathbf{T}^n(A_1, A_2)) \\
&\le M_{\max} \rho^{n - 1} H((B_1, B_2), (A_1, A_2)) \to 0
\end{align*}
as $n \to \infty$.
\end{proof}

To determine the Hausdorff dimension of the coupled attractors, we recall the following lemma.

\begin{lemma}[{\cite[Lemma 5.3]{Hutchinson1981}}]\label{lem:hutchinson_ball}
Let $\{V_i\}$ be a collection of disjoint open subsets of $\mathbb{R}^m$, and let $U$ be a ball of radius $r$. Suppose there exist constants $0 < c_1 \le c_2 < \infty$ such that every set $V_i$ contains a ball of radius $c_1 r$ and is contained in a ball of radius $c_2 r$. Then the number of closures $\overline{V_i}$ that intersect $U$ is at most $(1 + 2c_2)^m c_1^{-m}$.
\end{lemma}

\begin{definition}\label{def:cosc}
A CD-IFS $(\mathcal{F}, \tilde{\mathcal{F}}, \mathcal{G}, \tilde{\mathcal{G}})$ on $\mathbb{R}^m$ satisfies the \emph{Coupled Open Set Condition (COSC)} if there exist non-empty bounded open sets $U_1, U_2 \subset \mathbb{R}^m$ such that
\[
\bigcup_{i=1}^{N_1} f_i(U_2) \cup \bigcup_{k=1}^p g_k(U_1) \subseteq U_1 \quad \text{and} \quad \bigcup_{j=1}^{N_2} \tilde{f}_j(U_1) \cup \bigcup_{l=1}^q \tilde{g}_l(U_2) \subseteq U_2,
\]
where the sets $f_i(U_2)$ ($1 \le i \le N_1$) and $g_k(U_1)$ ($1 \le k \le p$) are all pairwise disjoint, and the sets $\tilde{f}_j(U_1)$ ($1 \le j \le N_2$) and $\tilde{g}_l(U_2)$ ($1 \le l \le q$) are all pairwise disjoint.
\end{definition}

Let $M_i$ and $\tilde{M}_j$ denote the linear parts of $f_i$ and $\tilde{f}_j$, respectively, and define the scaling moduli
\[
\alpha_i = |\det(M_i)|^{1/m} \quad (1 \le i \le N_1), \qquad \beta_j = |\det(\tilde{M}_j)|^{1/m} \quad (1 \le j \le N_2).
\]
Since $f_i \circ \tilde{f}_j$ is a similarity with ratio $c_{i,j}$, the determinant of its linear part satisfies $|\det(M_i \tilde{M}_j)| = c_{i,j}^m$, which gives
\begin{equation}\label{eq:c_factor}
c_{i,j} = |\det(M_i \tilde{M}_j)|^{1/m} = |\det(M_i)|^{1/m}|\det(\tilde{M}_j)|^{1/m} = \alpha_i \beta_j.
\end{equation}

\begin{theorem}\label{thm:hybrid_dimension}
Let $\mathcal{F} = \{f_1, \dots, f_{N_1}\}$, $\tilde{\mathcal{F}} = \{\tilde{f}_1, \dots, \tilde{f}_{N_2}\}$, $\mathcal{G} = \{g_1, \dots, g_p\}$, and $\tilde{\mathcal{G}} = \{\tilde{g}_1, \dots, \tilde{g}_q\}$ satisfy the hypotheses of Lemma~\ref{lem:sandwich} and the Coupled Open Set Condition, with coupled attractor $(A_1, A_2)$. Then the Hausdorff dimension $s = \dim_H(A_1) = \dim_H(A_2)$ is the unique real number in the region $\sum_{k=1}^p r_k^s < 1$, $\sum_{l=1}^q \tilde{r}_l^s < 1$ satisfying
\begin{equation}\label{eq:main_dim}
\left(1 - \sum_{k=1}^p r_k^s\right)\left(1 - \sum_{l=1}^q \tilde{r}_l^s\right) = \left(\sum_{i=1}^{N_1} \alpha_i^s\right)\left(\sum_{j=1}^{N_2} \beta_j^s\right).
\end{equation}
Moreover, $0 < \mathcal{H}^s(A_1) < \infty$ and $0 < \mathcal{H}^s(A_2) < \infty$.
\end{theorem}

\begin{proof}
Set $R_1(t) = \sum_{k=1}^p r_k^t$, $R_2(t) = \sum_{l=1}^q \tilde{r}_l^t$, $\Phi_1(t) = \sum_{i=1}^{N_1} \alpha_i^t$, and $\Phi_2(t) = \sum_{j=1}^{N_2} \beta_j^t$. If $p = q = 0$ and $N_1 = N_2 = 1$, then $A_1 = (f_1 \circ \tilde{f}_1)(A_1)$ and $A_2 = (\tilde{f}_1 \circ f_1)(A_2)$ are singletons, so $\dim_H(A_1) = \dim_H(A_2) = 0$ and $\mathcal{H}^0(A_1) = \mathcal{H}^0(A_2) = 1$, matching the unique real solution $s = 0$ to \eqref{eq:main_dim}. Assume now that $p + q \ge 1$ or $N_1 N_2 \ge 2$. Note that by \eqref{eq:c_factor}, $\Phi_1(t)\Phi_2(t) = \sum_{i,j} c_{i,j}^t$ is continuous and strictly decreasing on $\mathbb{R}$, with $\lim_{t \to \infty} \Phi_1(t)\Phi_2(t) = 0$.
\begin{itemize}
    \item If $p = q = 0$ and $N_1 N_2 \ge 2$, then \eqref{eq:main_dim} reduces to $\sum_{i,j} c_{i,j}^t = 1$. Since this sum equals $N_1 N_2 \ge 2$ at $t = 0$ and tends to $0$ as $t \to \infty$, it has a unique real solution $s > 0$.
    \item If $p + q \ge 1$, let $s_1 \ge 0$ solve $R_1(s_1) = 1$ if $p \ge 1$ (with $s_1 = 0$ if $p = 0$), let $s_2 \ge 0$ solve $R_2(s_2) = 1$ if $q \ge 1$ (with $s_2 = 0$ if $q = 0$), and set $t_0 = \max\{s_1, s_2\} \ge 0$. Because $R_1(0) = p$ and $R_2(0) = q$, the region $\{t \in \mathbb{R} : R_1(t) < 1 \text{ and } R_2(t) < 1\}$ is $(t_0, \infty)$. On $[t_0, \infty)$, the product $(1 - R_1(t))(1 - R_2(t))$ is continuous and strictly increasing from $0$ to $1$, while $\Phi_1(t)\Phi_2(t)$ decreases continuously from $\Phi_1(t_0)\Phi_2(t_0) > 0$ to $0$. By the Intermediate Value Theorem and monotonicity, there is a unique $s > t_0 \ge 0$ satisfying \eqref{eq:main_dim}, with $R_1(s) < 1$ and $R_2(s) < 1$.
\end{itemize}

In either case we have $s > 0$, $R_1(s) < 1$, and $R_2(s) < 1$. We prove $0 < \mathcal{H}^s(A_1) < \infty$; the argument for $A_2$ is symmetric. By \eqref{eq:main_dim}, we may define
\[
\gamma = \left(\frac{1 - R_1(s)}{\Phi_1(s)}\right)^{1/s} = \left(\frac{\Phi_2(s)}{1 - R_2(s)}\right)^{1/s} > 0,
\]
and assign a weight $w(\psi)$ to each symbol $\psi \in \mathcal{F} \cup \tilde{\mathcal{F}} \cup \mathcal{G} \cup \tilde{\mathcal{G}}$ by
\[
w(f_i) = \gamma \alpha_i, \quad w(\tilde{f}_j) = \gamma^{-1} \beta_j, \quad w(g_k) = r_k, \quad w(\tilde{g}_l) = \tilde{r}_l,
\]
with $w(\mathrm{Id}) = 1$. By the choice of $\gamma$, these weights satisfy
\begin{equation}\label{eq:weight_norm}
\sum_{i=1}^{N_1} w(f_i)^s + \sum_{k=1}^p w(g_k)^s = 1 \quad \text{and} \quad \sum_{j=1}^{N_2} w(\tilde{f}_j)^s + \sum_{l=1}^q w(\tilde{g}_l)^s = 1.
\end{equation}
Since $s > 0$ and every term in \eqref{eq:weight_norm} is positive, $0 < w(\psi) \le 1$ for every symbol $\psi$. Furthermore, for every pair $(i, j)$, the factor $\gamma$ cancels: $w(f_i)w(\tilde{f}_j) = (\gamma \alpha_i)(\gamma^{-1} \beta_j) = \alpha_i \beta_j = c_{i,j}$. By Lemma~\ref{lem:sandwich}(1), the contraction ratio of every block $\Phi = f_i \circ \tilde{g}_{l_1} \circ \dots \circ \tilde{g}_{l_d} \circ \tilde{f}_j$ therefore equals the product $w(f_i) w(\tilde{g}_{l_1}) \cdots w(\tilde{g}_{l_d}) w(\tilde{f}_j)$.

For any admissible word $\omega = (\psi_1, \dots, \psi_L) \in \mathcal{I}_1^L$, set $W(\omega) = \prod_{t=1}^L w(\psi_t)$, with $W(\emptyset) = 1$ for the empty word. Applying \eqref{eq:inplace_factor} from Lemma~\ref{lem:sandwich}(2), write $h_\omega = h^{\mathrm{head}} \circ h^{\mathrm{tail}} \circ h^{\mathrm{sim}}$, and let $W_{\mathrm{head}}$ and $W_{\mathrm{sim}}$ denote the weight products of the subwords corresponding to $h^{\mathrm{head}}$ and $h^{\mathrm{sim}}$ (equal to $1$ if the subword is empty). Then $h^{\mathrm{head}}$ and $h^{\mathrm{sim}}$ are similarities with ratios $W_{\mathrm{head}}$ and $W_{\mathrm{sim}}$, $h^{\mathrm{tail}} \in \mathcal{T} = \{\mathrm{Id}\} \cup \mathcal{F}$, and $W(\omega) = W_{\mathrm{head}} w(h^{\mathrm{tail}}) W_{\mathrm{sim}}$. Let $c_{\max}, r_{\max}, \rho$ be as in the proof of Theorem~\ref{thm:hybrid_coupled_existence}. If $\omega \in \mathcal{I}_1^L$ contains $m_{\mathrm{sim}}$ symbols from $\mathcal{G} \cup \tilde{\mathcal{G}}$ and $k = L - m_{\mathrm{sim}}$ symbols from $\mathcal{F} \cup \tilde{\mathcal{F}}$, then $m_{\mathrm{sim}} + 2\lfloor k/2 \rfloor \ge L - 1$ and $w(h^{\mathrm{tail}}) \le 1$, whence
\[
W(\omega) = W_{\mathrm{head}} w(h^{\mathrm{tail}}) W_{\mathrm{sim}} \le (r_{\max})^{m_{\mathrm{sim}}} (c_{\max})^{\lfloor k/2 \rfloor} \le \rho^{L - 1}.
\]
In particular, $\max_{\omega \in \mathcal{I}_1^L} W(\omega) \le \rho^{L - 1} \to 0$ as $L \to \infty$.

We first show that $\mathcal{H}^s(A_1) < \infty$. Since $\mathbf{T}^L(A_1, A_2) = (A_1, A_2)$, the sets $\{h_\omega(A_{\tau(\omega)}) : \omega \in \mathcal{I}_1^L\}$ cover $A_1$. Letting $D = \max\{\diam(A_1), \diam(A_2)\}$ and $\kappa_{\max} = \max_{\psi \in \mathcal{T}} \frac{\sigma_{\max}(\psi)}{w(\psi)} < \infty$ (where $\sigma_{\max}(\psi)$ and $\sigma_{\min}(\psi)$ denote the largest and smallest singular values of the linear part of $\psi$), we obtain
\begin{align*}
\diam(h_\omega(A_{\tau(\omega)})) &= W_{\mathrm{head}} \diam\big(h^{\mathrm{tail}}(h^{\mathrm{sim}}(A_{\tau(\omega)}))\big) \\
&\le W_{\mathrm{head}} \sigma_{\max}(h^{\mathrm{tail}}) W_{\mathrm{sim}} \diam(A_{\tau(\omega)}) \\
&\le \kappa_{\max} D \, W_{\mathrm{head}} w(h^{\mathrm{tail}}) W_{\mathrm{sim}} = \kappa_{\max} D \, W(\omega).
\end{align*}
Summing the $s$-th powers of the diameters over $\omega \in \mathcal{I}_1^L$ gives
\[
\sum_{\omega \in \mathcal{I}_1^L} (\diam(h_\omega(A_{\tau(\omega)})))^s \le (\kappa_{\max} D)^s \sum_{\omega \in \mathcal{I}_1^L} W(\omega)^s.
\]
By \eqref{eq:weight_norm} and induction on $L$, we have $\sum_{\omega \in \mathcal{I}_1^L} W(\omega)^s = 1$ for all $L \ge 1$. Since $(\kappa_{\max} D)^s$ is independent of $L$ and $\max_{\omega \in \mathcal{I}_1^L} \diam(h_\omega(A_{\tau(\omega)})) \le \kappa_{\max} D \, \rho^{L - 1} \to 0$ as $L \to \infty$, it follows that $\mathcal{H}^s(A_1) \le (\kappa_{\max} D)^s < \infty$.

We now prove that $\mathcal{H}^s(A_1) > 0$. Let $\mathcal{I}_1^* = \bigcup_{L=1}^\infty \mathcal{I}_1^L$, and let $\Sigma_1$ denote the space of infinite sequences $\sigma = (\psi_1, \psi_2, \dots)$ such that every finite prefix $\sigma|_L := (\psi_1, \dots, \psi_L)$ belongs to $\mathcal{I}_1^L$ ($L \ge 1$), equipped with the product topology. By \eqref{eq:weight_norm}, assigning mass $W(\omega)^s$ to each cylinder set $[\omega] = \{\sigma \in \Sigma_1 : \sigma|_L = \omega\}$ ($\omega \in \mathcal{I}_1^L$, $L \ge 1$) is consistent across levels and extends uniquely to a Borel probability measure $\tilde{\mu}$ on $\Sigma_1$. For every $\sigma \in \Sigma_1$, $\{h_{\sigma|_L}(A_{\tau(\sigma|_L)})\}_{L=1}^\infty$ is a nested sequence of non-empty compact sets whose diameters satisfy $\diam(h_{\sigma|_L}(A_{\tau(\sigma|_L)})) \le \kappa_{\max} D \, \rho^{L - 1} \to 0$ as $L \to \infty$; hence $\bigcap_{L=1}^\infty h_{\sigma|_L}(A_{\tau(\sigma|_L)})$ is a singleton, defining a continuous map $\pi: \Sigma_1 \to \mathbb{R}^m$ with $\pi(\Sigma_1) = A_1$. Let $\mu = \tilde{\mu} \circ \pi^{-1}$ be the pushforward Borel probability measure on $\mathbb{R}^m$ supported on $A_1$.

Let $U$ be an arbitrary closed ball of radius $r \in (0, 1)$. For $\omega = (\psi_1, \dots, \psi_L) \in \mathcal{I}_1^*$ of length $|\omega| = L$, write $\omega^- = (\psi_1, \dots, \psi_{L-1})$ (with $\omega^- = \emptyset$ when $L = 1$), and define the cut-set
\[
\mathcal{Q} = \{\omega \in \mathcal{I}_1^* : W(\omega) \le r < W(\omega^-)\}.
\]
Because $0 < w(\psi) \le 1$ for every symbol $\psi$, the sequence $\{W(\sigma|_L)\}_{L \ge 1}$ is non-increasing along every $\sigma \in \Sigma_1$, while $\max_{\omega \in \mathcal{I}_1^L} W(\omega) \le \rho^{L - 1} \to 0$ as $L \to \infty$. Since $W(\emptyset) = 1 > r$, every $\sigma \in \Sigma_1$ has a unique prefix in $\mathcal{Q}$, and any $\omega \in \mathcal{Q}$ with $|\omega| \ge 2$ satisfies $r < W(\omega^-) \le \rho^{|\omega| - 2}$; hence $\mathcal{Q}$ is finite and prefix-free, with $\Sigma_1 = \bigcup_{\omega \in \mathcal{Q}} [\omega]$ as a disjoint union. Setting $w_{\min} = \min_\psi w(\psi) > 0$, the definition of $\mathcal{Q}$ gives $w_{\min} r < W(\omega) \le r$ for all $\omega \in \mathcal{Q}$.

Let $U_1, U_2 \subset \mathbb{R}^m$ be non-empty bounded open sets satisfying Definition~\ref{def:cosc}, and choose $0 < a_1 \le a_2 < \infty$ such that each $U_i$ contains a ball of radius $a_1$ and is contained in a ball of radius $a_2$. Since $U_1$ and $U_2$ are bounded, $\overline{U_1}, \overline{U_2} \in \mathcal{K}(\mathbb{R}^m)$. By continuity, $T_i(\overline{U_1}, \overline{U_2}) \subseteq \overline{U_i}$ for $i \in \{1, 2\}$, and hence by induction $T_i^{(n)}(\overline{U_1}, \overline{U_2}) \subseteq \overline{U_i}$ for all $n \ge 1$; letting $n \to \infty$ via Theorem~\ref{thm:hybrid_coupled_existence} gives $A_1 \subseteq \overline{U_1}$ and $A_2 \subseteq \overline{U_2}$. Consider any two distinct words $\omega = (\omega_1, \dots, \omega_L)$ and $\nu = (\nu_1, \dots, \nu_{L'})$ in $\mathcal{Q}$. Since $\mathcal{Q}$ is prefix-free, there is a first index $t$ ($1 \le t \le \min\{L, L'\}$) at which $\omega_t \neq \nu_t$. Because $\omega$ and $\nu$ agree before position~$t$, both $\omega_t$ and $\nu_t$ belong to the same collection ($\mathcal{F} \cup \mathcal{G}$ or $\tilde{\mathcal{F}} \cup \tilde{\mathcal{G}}$). Repeated application of Definition~\ref{def:cosc} places the images of $U_{\tau(\omega)}$ and $U_{\tau(\nu)}$ under the remaining tails (taken to be $\mathrm{Id}$ if $t = L$ or $t = L'$) inside $U_{\tau(\omega_t)}$ and $U_{\tau(\nu_t)}$, respectively; hence their images under $\omega_t$ and $\nu_t$ are disjoint by the COSC, and applying the injective map $\omega_1 \circ \dots \circ \omega_{t-1}$ (taken to be $\mathrm{Id}$ if $t = 1$) preserves disjointness. Since each $h_\omega$ is an invertible affine map, $\{h_\omega(U_{\tau(\omega)}) : \omega \in \mathcal{Q}\}$ is a collection of pairwise disjoint open subsets of $U_1$, and for every $\omega \in \mathcal{Q}$ we have $\pi([\omega]) \subseteq h_\omega(A_{\tau(\omega)}) \subseteq h_\omega(\overline{U_{\tau(\omega)}}) = \overline{h_\omega(U_{\tau(\omega)})}$.

Set $\kappa_{\min} = \min_{\psi \in \mathcal{T}} \frac{\sigma_{\min}(\psi)}{w(\psi)} > 0$, $k_1 = \kappa_{\min} w_{\min} a_1$, and $k_2 = \kappa_{\max} a_2$, so that $0 < k_1 \le k_2 < \infty$. For each $\omega \in \mathcal{Q}$, writing $h_\omega = h^{\mathrm{head}} \circ h^{\mathrm{tail}} \circ h^{\mathrm{sim}}$ as in \eqref{eq:inplace_factor}, the maps $h^{\mathrm{sim}}$ and $h^{\mathrm{head}}$ are similarities with ratios $W_{\mathrm{sim}}$ and $W_{\mathrm{head}}$, while $h^{\mathrm{tail}}$ scales distances by at least $\sigma_{\min}(h^{\mathrm{tail}})$ and at most $\sigma_{\max}(h^{\mathrm{tail}})$. Since $U_{\tau(\omega)}$ contains a ball of radius $a_1$ and is contained in a ball of radius $a_2$, the image $h_\omega(U_{\tau(\omega)})$ contains a ball of radius
\[
W_{\mathrm{head}} \sigma_{\min}(h^{\mathrm{tail}}) W_{\mathrm{sim}} a_1 = \frac{\sigma_{\min}(h^{\mathrm{tail}})}{w(h^{\mathrm{tail}})} W(\omega) a_1 \ge (\kappa_{\min} w_{\min} a_1) r = k_1 r,
\]
and is contained in a ball of radius
\[
W_{\mathrm{head}} \sigma_{\max}(h^{\mathrm{tail}}) W_{\mathrm{sim}} a_2 = \frac{\sigma_{\max}(h^{\mathrm{tail}})}{w(h^{\mathrm{tail}})} W(\omega) a_2 \le (\kappa_{\max} a_2) r = k_2 r.
\]
Hence every $h_\omega(U_{\tau(\omega)})$ ($\omega \in \mathcal{Q}$) contains a ball of radius $k_1 r$ and is contained in a ball of radius $k_2 r$.

Setting $K = (1 + 2k_2)^m k_1^{-m}$, Lemma~\ref{lem:hutchinson_ball} shows that the set $\mathcal{Q}_U = \{\omega \in \mathcal{Q} : \overline{h_\omega(U_{\tau(\omega)})} \cap U \neq \emptyset\}$ has at most $K$ elements. Since $\Sigma_1 = \bigcup_{\omega \in \mathcal{Q}} [\omega]$ and $\pi([\omega]) \subseteq \overline{h_\omega(U_{\tau(\omega)})}$ for every $\omega \in \mathcal{Q}$, we have $\pi^{-1}(U) \subseteq \bigcup_{\omega \in \mathcal{Q}_U} [\omega]$, and hence
\[
\mu(U) = \tilde{\mu}(\pi^{-1}(U)) \le \sum_{\omega \in \mathcal{Q}_U} \tilde{\mu}([\omega]) = \sum_{\omega \in \mathcal{Q}_U} W(\omega)^s \le K r^s.
\]
In particular, for every $x \in \mathbb{R}^m$, applying this estimate to the closed ball of center $x$ and radius $r \in (0, 1)$ and letting $r \to 0^+$ yields $\mu(\{x\}) = 0$ (since $s > 0$). For any Borel set $E \subset \mathbb{R}^m$ with $0 < \diam(E) \le 1/2$, choosing $x_0 \in E$ places $E$ inside the closed ball of center $x_0$ and radius $\diam(E) < 1$, whence $\mu(E) \le K(\diam(E))^s$. Thus $\mu(E) \le K(\diam(E))^s$ holds for every Borel set $E \subset \mathbb{R}^m$ with $\diam(E) \le 1/2$. Since $\mu(A_1) = 1$, the Mass Distribution Principle \cite[Principle 4.2]{Falconer2014} gives $\mathcal{H}^s(A_1) \ge \mu(A_1)/K = 1/K > 0$. Combining $0 < \mathcal{H}^s(A_1) < \infty$ with the symmetric bounds for $A_2$ yields $\dim_H(A_1) = \dim_H(A_2) = s$.
\end{proof}

\begin{corollary}\label{cor:ahlfors_dims}
Under the hypotheses of Theorem~\ref{thm:hybrid_dimension}, the Hausdorff, packing, and box-counting dimensions of both coupled attractors coincide:
\[
\dim_H(A_1) = \dim_P(A_1) = \dim_B(A_1) = \dim_H(A_2) = \dim_P(A_2) = \dim_B(A_2) = s.
\]
Moreover, whenever $s > 0$, the pushforward measure $\mu$ supported on $A_1$ constructed in the proof of Theorem~\ref{thm:hybrid_dimension} (and symmetrically for $A_2$) is Ahlfors $s$-regular: there exist constants $0 < C_1 \le C_2 < \infty$ such that
\[
C_1 r^s \le \mu(B(x, r)) \le C_2 r^s \quad \text{for all } x \in A_1 \text{ and } 0 < r \le \diam(A_1),
\]
where $B(x, r)$ denotes the open ball of center $x$ and radius $r$.
\end{corollary}

\begin{proof}
If $s = 0$, both attractors are singletons and all dimensions equal $0$. Assume $s > 0$. By Theorem~\ref{thm:hybrid_dimension}, $\mathcal{H}^s(A_1) > 0$, so $\diam(A_1) > 0$. Using the constants $K, D, \kappa_{\max}, w_{\min}$ from the proof of Theorem~\ref{thm:hybrid_dimension}, set $C_2 = \max\{K, 1\}$. For any $x \in A_1$, letting $\overline{B}(x, r)$ denote the closed ball of center $x$ and radius $r$, we have $\mu(B(x, r)) \le \mu(\overline{B}(x, r)) \le K r^s \le C_2 r^s$ when $0 < r < 1$, and $\mu(B(x, r)) \le 1 \le r^s \le C_2 r^s$ when $r \ge 1$; thus $\mu(B(x, r)) \le C_2 r^s$ for all $r > 0$.

For the lower bound, set $C_1 = \left(\frac{w_{\min}}{2\kappa_{\max} D}\right)^s > 0$, and fix $x \in A_1$ and $0 < r \le \diam(A_1)$. Since $\pi(\Sigma_1) = A_1$, choose $\sigma \in \Sigma_1$ such that $\pi(\sigma) = x$, and set $\eta = \frac{r}{2\kappa_{\max} D}$. Because $\kappa_{\max} \ge \sigma_{\max}(\mathrm{Id})/w(\mathrm{Id}) = 1$ and $D \ge \diam(A_1) \ge r$, we have $0 < \eta \le 1/2 < 1$. Hence there is a unique prefix $\omega$ of $\sigma$ satisfying $W(\omega) \le \eta < W(\omega^-)$, which gives $W(\omega) > w_{\min}\eta = \frac{w_{\min}}{2\kappa_{\max} D} r$. Since $x \in h_\omega(A_{\tau(\omega)})$ and $\diam(h_\omega(A_{\tau(\omega)})) \le \kappa_{\max} D W(\omega) \le r/2 < r$, we have $\pi([\omega]) \subseteq h_\omega(A_{\tau(\omega)}) \subseteq B(x, r)$, so $[\omega] \subseteq \pi^{-1}(B(x, r))$ and
\[
\mu(B(x, r)) = \tilde{\mu}(\pi^{-1}(B(x, r))) \ge \tilde{\mu}([\omega]) = W(\omega)^s > \left(\frac{w_{\min}}{2\kappa_{\max} D}\right)^s r^s = C_1 r^s.
\]
In particular, $0 < C_1 \le C_2 < \infty$, which establishes Ahlfors $s$-regularity.

Finally, for any $0 < r \le \diam(A_1)$, since $A_1$ is compact, we can choose a finite maximal pairwise disjoint collection $\{B(x_j, r/2)\}_{j=1}^{N_r}$ of open balls with centers $x_j \in A_1$. Applying the lower bound at radius $r/2$ gives
\[
1 = \mu(\mathbb{R}^m) \ge \sum_{j=1}^{N_r} \mu(B(x_j, r/2)) \ge N_r C_1 \left(\frac{r}{2}\right)^s,
\]
whence $N_r \le C_1^{-1} 2^s r^{-s}$. By maximality, for every $y \in A_1$ there exists $j$ such that $B(y, r/2) \cap B(x_j, r/2) \neq \emptyset$, so $\|y - x_j\| < r$. Thus the balls $\{B(x_j, r)\}_{j=1}^{N_r}$ of diameter $2r$ cover $A_1$, yielding $\overline{\dim}_B(A_1) \le \limsup_{r \to 0^+} \frac{\ln N_r}{-\ln(2r)} \le s$. Combining this with the standard inequalities $s = \dim_H(A_1) \le \underline{\dim}_B(A_1) \le \overline{\dim}_B(A_1) \le s$ and $s = \dim_H(A_1) \le \dim_P(A_1) \le \overline{\dim}_B(A_1) \le s$ proves the equality of all dimensions for $A_1$, and symmetrically for $A_2$.
\end{proof}

We now apply Theorem~\ref{thm:hybrid_dimension} to explicit coupled self-affine systems and verify the COSC in each case.

\begin{figure}[htbp]
    \centering
    \includegraphics[width=\textwidth]{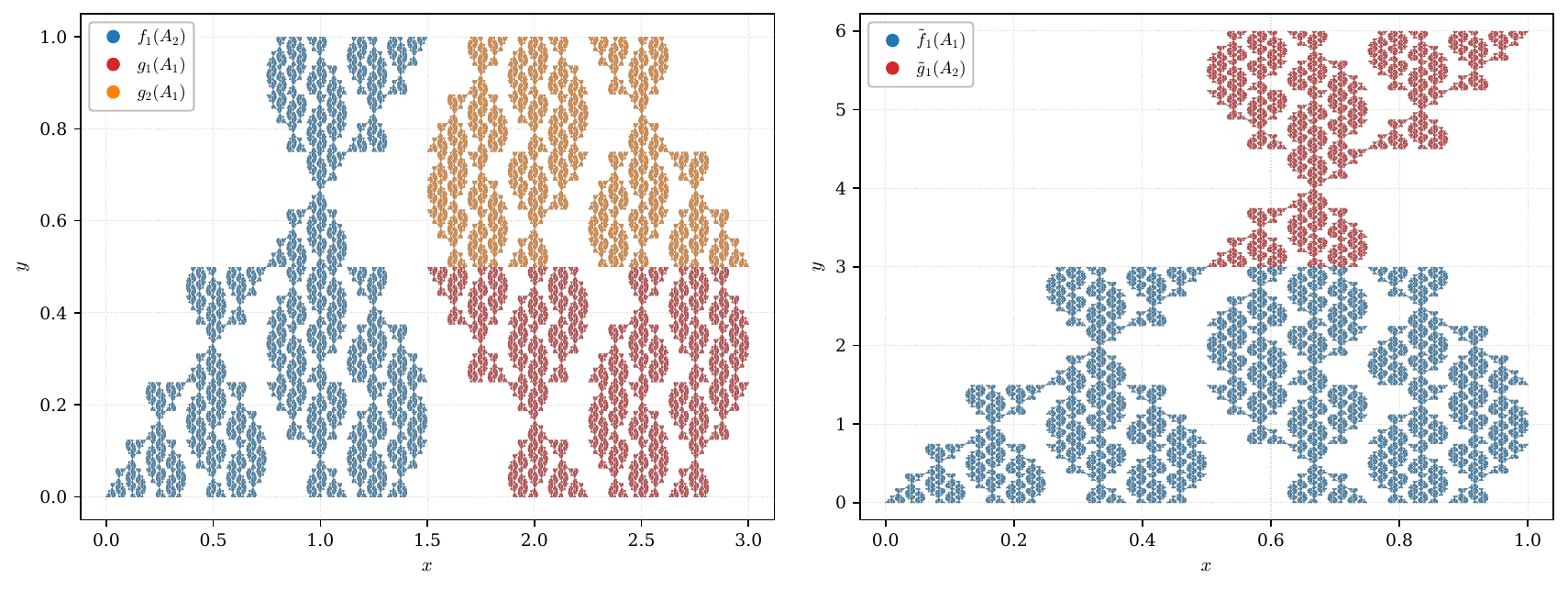}
    \caption{Coupled attractors $A_1$ (left) and $A_2$ (right) of Example~\ref{ex:area_preserving}, colored by their first-level images ($\dim_H(A_1) = \dim_H(A_2) = \log_2(2 + \sqrt{2}) \approx 1.7716$).}
    \label{fig:area_preserving}
\end{figure}

\begin{example}\label{ex:area_preserving}
Consider the coupled system on $\mathbb{R}^2$ generated by
\begin{gather*}
f_1(x, y) = \left(\frac{3}{2}x, \, \frac{y}{6}\right), \qquad \tilde{f}_1(x, y) = \left(\frac{x}{3}, \, 3y\right), \\
g_1(x, y) = \left(\frac{x}{2} + \frac{3}{2}, \, -\frac{y}{2} + \frac{1}{2}\right), \qquad g_2(x, y) = \left(-\frac{x}{2} + 3, \, \frac{y}{2} + \frac{1}{2}\right), \\
\tilde{g}_1(x, y) = \left(-\frac{x}{2} + 1, \, -\frac{y}{2} + 6\right).
\end{gather*}
Here $\tilde{f}_1$ is area-preserving with $\det(\tilde{M}_1) = 1$ ($\beta_1 = 1$), whereas $f_1$ has $\det(M_1) = 1/4$ ($\alpha_1 = 1/2$). Since $\Lip(f_1^n) = (3/2)^n > 1$ and $\Lip(\tilde{f}_1^n) = 3^n > 1$ for all $n \ge 1$, while both compositions $f_1 \circ \tilde{f}_1$ and $\tilde{f}_1 \circ f_1$ are similarities with ratio $c = \alpha_1 \beta_1 = 1/2$, the pair $(f_1, \tilde{f}_1)$ is a proper similarity contraction divisor pair. Moreover, since all linear parts are diagonal and therefore commute, the similarities $g_1, g_2$ ($r_1 = r_2 = 1/2$) are universally aligned with $\tilde{\mathcal{F}}$ and $\tilde{g}_1$ ($\tilde{r}_1 = 1/2$) is universally aligned with $\mathcal{F}$.

For the open rectangles $U_1 = (0, 3) \times (0, 1)$ and $U_2 = (0, 1) \times (0, 6)$, the images in State~$1$ are
\[
f_1(U_2) = \left(0, \frac{3}{2}\right) \times (0, 1), \quad g_1(U_1) = \left(\frac{3}{2}, 3\right) \times \left(0, \frac{1}{2}\right), \quad g_2(U_1) = \left(\frac{3}{2}, 3\right) \times \left(\frac{1}{2}, 1\right),
\]
which are pairwise disjoint subsets of $U_1$, while in State~$2$ the images
\[
\tilde{f}_1(U_1) = (0, 1) \times (0, 3), \qquad \tilde{g}_1(U_2) = \left(\frac{1}{2}, 1\right) \times (3, 6)
\]
are disjoint subsets of $U_2$, separated by the line $y = 3$. Thus the COSC holds. By Theorem~\ref{thm:hybrid_dimension} and Corollary~\ref{cor:ahlfors_dims}, the common Hausdorff, packing, and box-counting dimension $s$ is the unique solution in the region $2(1/2)^s < 1$, $(1/2)^s < 1$ of
\[
\left(1 - 2\left(\frac{1}{2}\right)^s\right)\left(1 - \left(\frac{1}{2}\right)^s\right) = \left(\frac{1}{2}\right)^s \cdot 1^s.
\]
Setting $u = (1/2)^s \in (0, 1/2)$ gives $2u^2 - 4u + 1 = 0$, whose unique root in $(0, 1/2)$ is $u = 1 - \frac{1}{\sqrt{2}}$. Hence $s = \log_2(2 + \sqrt{2}) \approx 1.7716$. The coupled attractors $(A_1, A_2)$ are shown in Figure~\ref{fig:area_preserving}.
\end{example}

\begin{figure}[htbp]
    \centering
    \includegraphics[width=\textwidth]{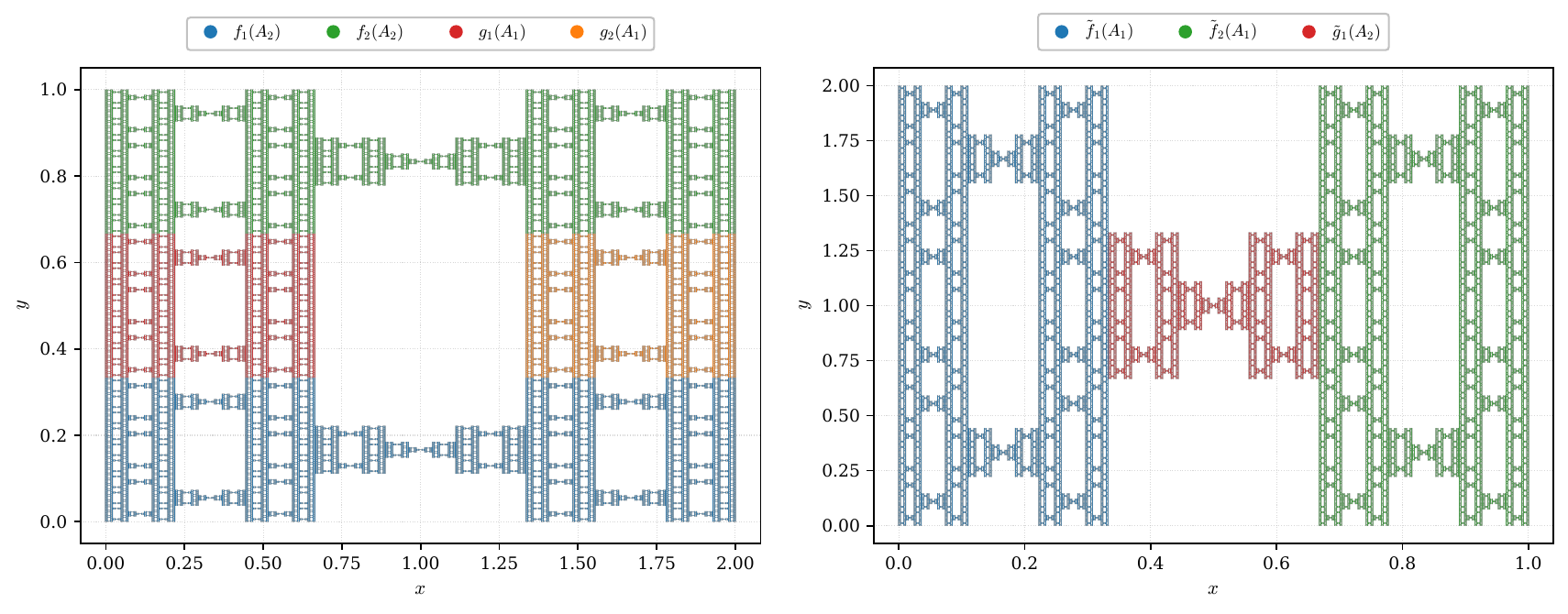}
    \caption{Coupled attractors $A_1$ (left) and $A_2$ (right) of Example~\ref{ex:hybrid_calc}, colored by their first-level images ($\dim_H(A_1) = \dim_H(A_2) = \log_3\big(\frac{7 + \sqrt{41}}{2}\big) \approx 1.7316$).}
    \label{fig:woven_carpet}
\end{figure}

\begin{example}\label{ex:hybrid_calc}
Consider the coupled system on $\mathbb{R}^2$ with $N_1 = 2$, $N_2 = 2$, $p = 2$, and $q = 1$, generated by
\begin{gather*}
f_1(x, y) = \left(2x, \, \frac{y}{6}\right), \qquad f_2(x, y) = \left(-2x + 2, \, -\frac{y}{6} + 1\right), \\
\tilde{f}_1(x, y) = \left(\frac{x}{6}, \, 2y\right), \qquad \tilde{f}_2(x, y) = \left(-\frac{x}{6} + 1, \, -2y + 2\right), \\
g_1(x, y) = \left(\frac{x}{3}, \, \frac{y}{3} + \frac{1}{3}\right), \qquad g_2(x, y) = \left(\frac{x}{3} + \frac{4}{3}, \, -\frac{y}{3} + \frac{2}{3}\right), \\
\tilde{g}_1(x, y) = \left(-\frac{x}{3} + \frac{2}{3}, \, \frac{y}{3} + \frac{2}{3}\right).
\end{gather*}
Since $\Lip(f_i^n) = \Lip(\tilde{f}_j^n) = 2^n > 1$ for all $n \ge 1$, while $f_i \circ \tilde{f}_j$ and $\tilde{f}_j \circ f_i$ are similarities with ratio $c_{i,j} = 1/3$, every pair $(f_i, \tilde{f}_j)$ ($i, j \in \{1, 2\}$) is a proper similarity contraction divisor pair with $\alpha_i = \beta_j = 1/\sqrt{3}$. Moreover, since all linear parts are diagonal and therefore commute, the similarities $g_1, g_2$ ($r_1 = r_2 = 1/3$) are universally aligned with $\tilde{\mathcal{F}}$ and $\tilde{g}_1$ ($\tilde{r}_1 = 1/3$) is universally aligned with $\mathcal{F}$.

Taking $U_1 = (0, 2) \times (0, 1)$ and $U_2 = (0, 1) \times (0, 2)$, we compute
\begin{gather*}
f_1(U_2) = (0, 2) \times \left(0, \frac{1}{3}\right), \qquad f_2(U_2) = (0, 2) \times \left(\frac{2}{3}, 1\right), \\
g_1(U_1) = \left(0, \frac{2}{3}\right) \times \left(\frac{1}{3}, \frac{2}{3}\right), \qquad g_2(U_1) = \left(\frac{4}{3}, 2\right) \times \left(\frac{1}{3}, \frac{2}{3}\right),
\end{gather*}
which are pairwise disjoint subsets of $U_1$. Similarly,
\[
\tilde{f}_1(U_1) = \left(0, \frac{1}{3}\right) \times (0, 2), \quad \tilde{f}_2(U_1) = \left(\frac{2}{3}, 1\right) \times (0, 2), \quad \tilde{g}_1(U_2) = \left(\frac{1}{3}, \frac{2}{3}\right) \times \left(\frac{2}{3}, \frac{4}{3}\right)
\]
are pairwise disjoint subsets of $U_2$. Thus the COSC holds, and Theorem~\ref{thm:hybrid_dimension} and Corollary~\ref{cor:ahlfors_dims} give the common Hausdorff, packing, and box-counting dimension $s$ as the unique solution in the region $2(1/3)^s < 1$, $(1/3)^s < 1$ of
\[
\left(1 - 2\left(\frac{1}{3}\right)^s\right)\left(1 - \left(\frac{1}{3}\right)^s\right) = \left(2\left(\frac{1}{\sqrt{3}}\right)^s\right)\left(2\left(\frac{1}{\sqrt{3}}\right)^s\right) = 4\left(\frac{1}{3}\right)^s.
\]
Setting $u = (1/3)^s \in (0, 1/2)$ yields $2u^2 - 7u + 1 = 0$, with unique root $u = \frac{7 - \sqrt{41}}{4} \in (0, 1/2)$. Therefore, $s = \log_3\left(\frac{7 + \sqrt{41}}{2}\right) \approx 1.7316$. The coupled attractors $(A_1, A_2)$ are shown in Figure~\ref{fig:woven_carpet}.
\end{example}

\begin{figure}[htbp]
    \centering
    \includegraphics[width=\textwidth]{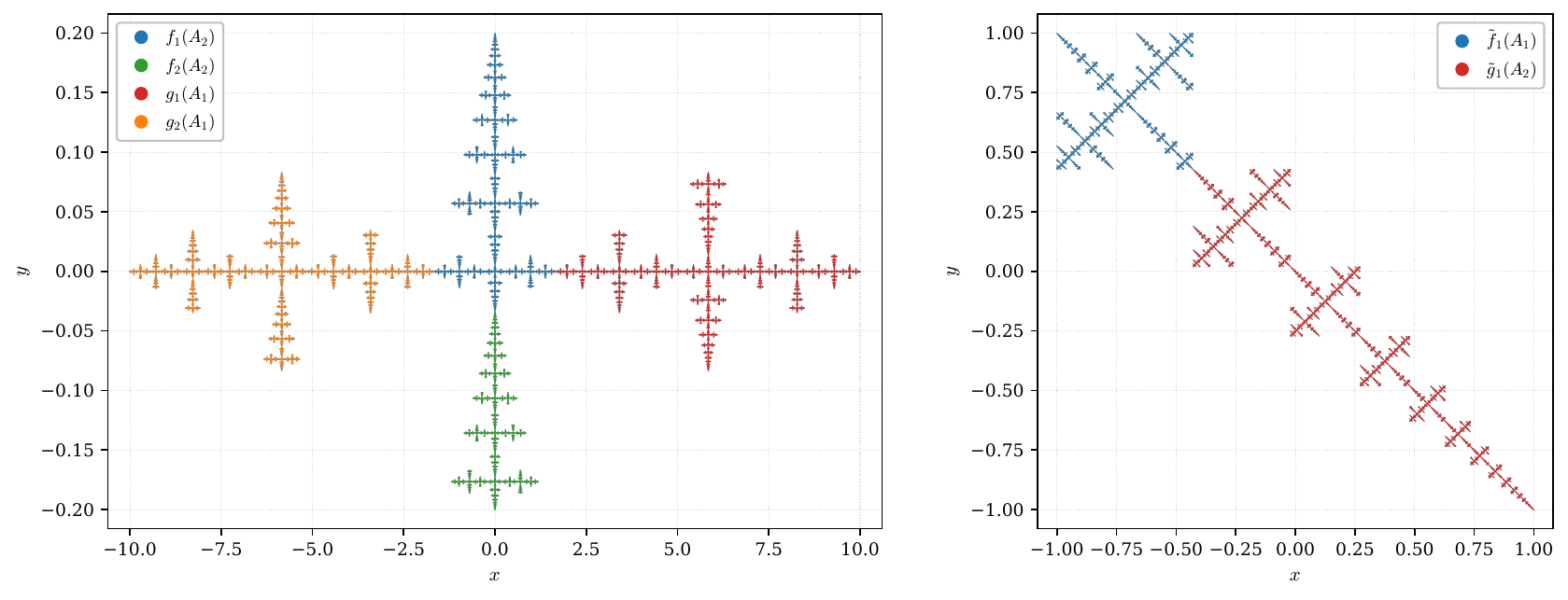}
    \caption{Coupled attractors $A_1$ (left) and $A_2$ (right) of Example~\ref{ex:noncommuting}, colored by their first-level images ($\dim_H(A_1) = \dim_H(A_2) \approx 1.3579$).}
    \label{fig:noncommuting}
\end{figure}

\begin{example}\label{ex:noncommuting}
Consider the coupled system on $\mathbb{R}^2$ with $N_1 = 2$, $N_2 = 1$, $p = 2$, and $q = 1$, generated by
\begin{gather*}
f_1(x, y) = \left(-\frac{35}{12}(x + y), \, \frac{7}{120}(x - y) + \frac{1}{12}\right), \\
f_2(x, y) = \left(-\frac{25}{12}(x + y), \, \frac{x - y}{24} - \frac{7}{60}\right), \\
\tilde{f}_1(x, y) = \left(-\frac{x}{35} - \frac{10}{7}y - \frac{5}{7}, \, -\frac{x}{35} + \frac{10}{7}y + \frac{5}{7}\right), \\
g_1(x, y) = \left(\frac{5}{12}x + \frac{35}{6}, \, -\frac{5}{12}y\right), \qquad g_2(x, y) = \left(\frac{5}{12}x - \frac{35}{6}, \, \frac{5}{12}y\right), \\
\tilde{g}_1(x, y) = \left(-\frac{5}{7}y + \frac{2}{7}, \, -\frac{5}{7}x - \frac{2}{7}\right).
\end{gather*}
The linear parts $M_1, M_2$, and $\tilde{M}_1$ of $f_1, f_2$, and $\tilde{f}_1$ have spectral radii $7(51 + \sqrt{2201})/240 > 1$, $(51 + \sqrt{2201})/48 > 1$, and $(49 + \sqrt{2801})/70 > 1$, respectively, so no positive iterate of $f_1, f_2$, or $\tilde{f}_1$ is a contraction. Composing the maps gives
\begin{gather*}
(f_1 \circ \tilde{f}_1)(x, y) = \left(\frac{x}{6}, \, -\frac{y}{6}\right), \quad (\tilde{f}_1 \circ f_1)(x, y) = \left(\frac{y - 5}{6}, \, \frac{x + 5}{6}\right), \\
(f_2 \circ \tilde{f}_1)(x, y) = \left(\frac{5x}{42}, \, -\frac{5y}{42} - \frac{37}{210}\right), \quad (\tilde{f}_1 \circ f_2)(x, y) = \left(\frac{5y - 23}{42}, \, \frac{5x + 23}{42}\right).
\end{gather*}
Thus $(f_1, \tilde{f}_1)$ and $(f_2, \tilde{f}_1)$ are proper, non-commuting similarity contraction divisor pairs with distinct ratios $c_{1,1} = 1/6$ and $c_{2,1} = 5/42$ ($\alpha_1 = 7/12$, $\alpha_2 = 5/12$, $\beta_1 = 2/7$). Since $\tilde{M}_1^\top \tilde{M}_1$ is diagonal, it commutes with the diagonal linear parts of $g_1$ and $g_2$ ($r_1 = r_2 = 5/12$); similarly, $M_1^\top M_1$ and $M_2^\top M_2$ are symmetric with equal diagonal entries, so they commute with the linear part of $\tilde{g}_1$ ($\tilde{r}_1 = 5/7$). Hence $\mathcal{G}$ is universally aligned with $\tilde{\mathcal{F}}$ and $\tilde{\mathcal{G}}$ is universally aligned with $\mathcal{F}$.

Consider the open rhombus $U_1 = \{(x, y) \in \mathbb{R}^2 : |x|/10 + 5|y| < 1\}$ and the open square $U_2 = (-1, 1) \times (-1, 1)$. In State~$1$, setting $(u, v) = f_1(x, y)$ gives $u/10 = -\frac{7}{24}(x + y)$ and $5v - \frac{5}{12} = \frac{7}{24}(x - y)$, so the identity $|a + b| + |a - b| = 2\max(|a|, |b|)$ yields $|u|/10 + |5v - \frac{5}{12}| = \frac{7}{12}\max(|x|, |y|)$. Arguing similarly for $f_2, g_1$, and $g_2$, we obtain
\begin{align*}
f_1(U_2) &= \left\{(x, y) : \frac{|x|}{10} + \left|5y - \frac{5}{12}\right| < \frac{7}{12}\right\}, \\
f_2(U_2) &= \left\{(x, y) : \frac{|x|}{10} + \left|5y + \frac{7}{12}\right| < \frac{5}{12}\right\}, \\
g_1(U_1) &= \left\{(x, y) : \left|\frac{x}{10} - \frac{7}{12}\right| + 5|y| < \frac{5}{12}\right\}, \\
g_2(U_1) &= \left\{(x, y) : \left|\frac{x}{10} + \frac{7}{12}\right| + 5|y| < \frac{5}{12}\right\}.
\end{align*}
By the triangle inequality, all four images lie inside $U_1$. Moreover,
\begin{gather*}
f_1(U_2) \subset \left\{5y + \frac{1}{6} > \frac{|x|}{10}\right\}, \qquad
f_2(U_2) \subset \left\{5y + \frac{1}{6} < -\frac{|x|}{10}\right\}, \\
g_1(U_1) \subset \left\{\frac{x}{10} > \left|5y + \frac{1}{6}\right|\right\}, \qquad
g_2(U_1) \subset \left\{\frac{x}{10} < -\left|5y + \frac{1}{6}\right|\right\},
\end{gather*}
so the four images are pairwise disjoint.

In State~$2$, setting $(u, v) = \tilde{f}_1(x, y)$ gives $-x/10 = \frac{7}{4}((u + 5/7) + (v - 5/7))$ and $5y = \frac{7}{4}((v - 5/7) - (u + 5/7))$, whence $|x|/10 + 5|y| = \frac{7}{2}\max(|u + 5/7|, |v - 5/7|)$. Therefore,
\[
\tilde{f}_1(U_1) = \left(-1, -\frac{3}{7}\right) \times \left(\frac{3}{7}, 1\right), \qquad \tilde{g}_1(U_2) = \left(-\frac{3}{7}, 1\right) \times \left(-1, \frac{3}{7}\right),
\]
which are disjoint open subsets of $U_2$. Thus the COSC holds.

By Theorem~\ref{thm:hybrid_dimension} and Corollary~\ref{cor:ahlfors_dims}, the common Hausdorff, packing, and box-counting dimension $s$ is the unique solution in the region $2(5/12)^s < 1$, $(5/7)^s < 1$ of
\[
\left(1 - 2\left(\frac{5}{12}\right)^s\right)\left(1 - \left(\frac{5}{7}\right)^s\right) = \left(\frac{1}{6}\right)^s + \left(\frac{5}{42}\right)^s,
\]
yielding $s \approx 1.3579$. The coupled attractors $(A_1, A_2)$ are shown in Figure~\ref{fig:noncommuting}.
\end{example}

\begin{example}\label{ex:non_aligned}
To show that universal alignment cannot be omitted in Theorem~\ref{thm:hybrid_dimension}, consider the coupled system on $\mathbb{R}^2$ with
\[
f_1(x, y) = \left(2x, \, \frac{y}{50}\right), \qquad \tilde{f}_1(x, y) = \left(\frac{x}{50}, \, 2y\right),
\]
and the $90^\circ$-rotated similarities of ratio $r_j = \tilde{r}_i = 1/4$:
\[
g_j(x, y) = \left(-\frac{y}{4} + 10\varepsilon_j, \, \frac{x}{4} + \varepsilon_j\right), \qquad \tilde{g}_i(x, y) = \left(-\frac{y}{4} + \varepsilon_i, \, \frac{x}{4} + 10\varepsilon_i\right),
\]
where $i, j \in \{1, 2\}$, $\varepsilon_1 = 1$, and $\varepsilon_2 = -1$. Here $(f_1, \tilde{f}_1)$ is a proper similarity contraction divisor pair with $\alpha_1 = \beta_1 = 1/5$ and $c_{1,1} = 1/25$ (as $\Lip(f_1^n) = \Lip(\tilde{f}_1^n) = 2^n > 1$ for all $n \ge 1$).

Although the $90^\circ$ rotations violate universal alignment (so Theorem~\ref{thm:hybrid_coupled_existence} does not apply directly), every admissible composition of length~$2$ is a contraction in the Euclidean metric: $\Lip(f_1 \circ \tilde{f}_1) = \Lip(\tilde{f}_1 \circ f_1) = 1/25$, $\Lip(g_j \circ g_k) = \Lip(\tilde{g}_i \circ \tilde{g}_k) = 1/16$, and each mixed composition ($f_1 \circ \tilde{g}_i$, $g_j \circ f_1$, $\tilde{f}_1 \circ g_j$, $\tilde{g}_i \circ \tilde{f}_1$) has Lipschitz constant at most $2 \cdot \frac{1}{4} = \frac{1}{2}$. Hence $\mathbf{T}^2$ is a contraction with factor $1/2$ in the product Hausdorff metric $H$, so it admits a unique fixed point $(A_1, A_2)$. Since $\mathbf{T}^2(\mathbf{T}(A_1, A_2)) = \mathbf{T}(\mathbf{T}^2(A_1, A_2)) = \mathbf{T}(A_1, A_2)$, uniqueness implies $\mathbf{T}(A_1, A_2) = (A_1, A_2)$, yielding a unique compact coupled attractor $(A_1, A_2)$.

The system satisfies the COSC with $U_1 = (-12, 12) \times (-4, 4)$ and $U_2 = (-4, 4) \times (-12, 12)$, since $f_1(U_2) = (-8, 8) \times (-6/25, 6/25)$, $g_1(U_1) = (9, 11) \times (-2, 4)$, and $g_2(U_1) = (-11, -9) \times (-4, 2)$ are subsets of $U_1$ with pairwise disjoint $x$-projections, and similarly $\tilde{f}_1(U_1) = (-6/25, 6/25) \times (-8, 8)$, $\tilde{g}_1(U_2) = (-2, 4) \times (9, 11)$, and $\tilde{g}_2(U_2) = (-4, 2) \times (-11, -9)$ are subsets of $U_2$ with pairwise disjoint $y$-projections. If \eqref{eq:main_dim} were applied despite the failure of universal alignment, it would give $(1 - 2(1/4)^s)^2 = (1/25)^s$ in the region $2(1/4)^s < 1$, or $(1/5)^s + 2(1/4)^s = 1$, which yields $s \approx 0.7542 < 1$.

However, evaluating the four admissible composite maps $\Phi_{i,j} = f_1 \circ \tilde{g}_i \circ \tilde{f}_1 \circ g_j$ ($i, j \in \{1, 2\}$) gives
\[
\Phi_{i,j}(x, y) = \left(-\frac{x}{4} + 2\varepsilon_i - \varepsilon_j, \, -\frac{y}{40000} + \frac{\varepsilon_i}{5} + \frac{\varepsilon_j}{1000}\right).
\]
The first coordinate of $\Phi_{i,j}$ depends only on $x$, and as $(\varepsilon_i, \varepsilon_j)$ ranges over $\{-1, 1\}^2$, the translation term $2\varepsilon_i - \varepsilon_j$ takes all four values in $\{-3, -1, 1, 3\}$. Since $A_1 \supseteq \bigcup_{i,j} \Phi_{i,j}(A_1)$, the projection $P_x(A_1)$ onto the $x$-axis satisfies $P_x(A_1) \supseteq \bigcup_{b \in \{-3,-1,1,3\}} \psi_b(P_x(A_1))$ for $\psi_b(x) = -x/4 + b$. Because $P_x(A_1)$ is non-empty and compact, iterating this inclusion and passing to the limit shows that $P_x(A_1)$ contains the attractor of $\{\psi_b\}$, which is the interval $[-4, 4]$. Therefore $\dim_H(A_1) \ge \dim_H(P_x(A_1)) = 1 > 0.7542$.
\end{example}

\section*{Acknowledgements}
This research was financially supported by the University Grants Commission (UGC), Government of India, through the award of the Junior Research Fellowship (JRF) (Ref.\ No.\ 1090 (CSIR-UGC NET JUNE 2019)) awarded to the first author. The first and second authors acknowledge the institutional support provided by Rajagiri School of Engineering and Technology, Kerala, India, and the third author acknowledges the institutional support provided by Muthoot Institute of Technology and Science, Kerala, India.

\bibliographystyle{amsplain}
\bibliography{references}

\end{document}